\input ./style/arxiv-general.cfg
\documentclass[aop,MSNbibl,seceqn,dvips]{arximspdf}
\makeatletter
   \@ifpackageloaded{graphicx}{}{\usepackage{graphicx}}
\makeatother
\doi{10.1214/14-AOP941} %kopijuoti is PTS
\volume{43}
\issue{5}
\pubyear{2015}
\firstpage{2481}
\lastpage{2510}
\docsubty{FLA}

\makeatletter
\newtheorem{teo}{Theorem}
\newproclaim{rem}[teo]{Remark}
\newtheorem{cor}[teo]{Corollary}
\newtheorem{prop}[teo]{Proposition}
\makeatother

\begin{document}
\begin{frontmatter}

%\dochead{}
\title{Embedding laws in diffusions by functions of time}
\runtitle{Embedding laws in diffusions by functions of time}

\begin{aug}
\author[A]{\fnms{A.~M.~G.}~\snm{Cox}\ead[label=e1]{a.m.g.cox@bath.ac.uk}}
\and
\author[B]{\fnms{G.}~\snm{Peskir}\corref{}\ead[label=e2]{goran@maths.man.ac.uk}}
\runauthor{A. M. G. Cox and G. Peskir}
\affiliation{University of Bath and The University of Manchester}
%\dedicated{}
\address[A]{Department of Mathematical Sciences\\
University of Bath\\
Claverton Down\\
Bath BA2 7AY\\
United Kingdom\\
\printead{e1}} %adresu isvedimo komanda gale!
\address[B]{School of Mathematics\\
The University of Manchester\\
Oxford Road\\
Manchester M13 9PL\\
United Kingdom\\
\printead{e2}}
\end{aug}

% HISTORY:
\received{\smonth{1} \syear{2013}}
\revised{\smonth{5} \syear{2014}}
%\accepted{\smonth{} \syear{}}

% ABSTRACT
%
\begin{abstract}
We present a constructive probabilistic proof of the fact that if
$B=(B_t)_{t \ge0}$ is standard Brownian motion started at $0$, and
$\mu$ is a given probability measure on $\mathbb{R}$ such that
$\mu(\{0\})=0$, then there exists a unique left-continuous increasing
function $b\dvtx (0,\infty) \rightarrow\mathbb{R}\cup\{+\infty\}$ and a
unique \mbox{left-}continuous decreasing function $c\dvtx (0,\infty)
\rightarrow\mathbb{R}\cup\{-\infty\}$ such that $B$ stopped at
$\tau_{b,c} = \inf\{ t>0 \vert B_t \ge b(t)$ or $B_t \le c(t) \}$
has the law $\mu$. The method of
proof relies upon weak convergence arguments arising from Helly's
selection theorem and makes use of the L\'evy metric which appears
to be novel in the context of embedding theorems. We show that
$\tau_{b,c}$ is minimal in the sense of Monroe so that the stopped
process $B^{\tau_{b,c}} = (B_{t \wedge\tau_{b,c}})_{t \ge0}$
satisfies natural uniform integrability conditions expressed in
terms of $\mu$. We also show that $\tau_{b,c}$ has the smallest
truncated expectation among all stopping times that embed $\mu$ into~$B$. The main results extend from standard Brownian motion to all
recurrent diffusion processes on the real line.
\end{abstract}

% KEYWORDS
% Pirmas kwd is didziosios raides
%
\end{frontmatter}

%s1 #&#
\section{Introduction}\label{sec1}
%%%%%%%%%%%%%%%%%%%%%%

A classic problem in modern probability theory is to find a stopping
time $\tau$ of a standard Brownian motion $B$ started at zero such
that $B$ stopped at $\tau$ has a given law $\mu$. The existence of a
randomised stopping time $\tau$ for centred laws $\mu$ was first
derived by Skorokhod \cite{Sk}, and the problem is often referred to
as the \emph{Skorokhod embedding problem}. A few years later Dubins
\cite{Du} proved the existence of a non-randomised stopping time
$\tau$ of $B$ that also holds for more general laws $\mu$. Many
other solutions have been found in subsequent years and we refer to
the survey article by Ob{\l}\'oj \cite{Ob} for a comprehensive
discussion (see also \cite{Ho} for financial applications and
\cite{KK} for connections to the Cantelli conjecture).

Solutions relevant to the present paper are those found by Root
\cite{Ro} in the setting of $B$ and Rost \cite{Ro-1} in the setting
of more general Markov processes and initial laws. Root \cite{Ro}
showed that $\tau$ can be realised as the first entry time to a
\emph{barrier}, and Rost \cite{Ro-1} showed that $\tau$ can be
characterised in terms of a \emph{filling scheme} dating back to
Chacon and Ornstein \cite{CO} within ergodic theory (see also
\cite{Di} for a closely related construction). Subsequently Chacon
\cite{Chh} showed that a stopping time arising from the filling
scheme coincides with the first entry time to a \emph{reversed
barrier} under some conditions. The proof of Root \cite{Ro} relies
upon a continuous mapping theorem and compactness of barriers in a
uniform distance topology. The methods of Rost \cite{Ro-1} and
Chacon \cite{Chh} rely on potential theory of general Markov
processes. Uniqueness of barriers was studied by Loynes \cite{Lo}.
He described barriers by functions of space. Reversed barriers can
also be described by functions of time. Based on this fact McConnell
\cite{Mc} developed an analytic free-boundary approach relying upon
potential theoretic considerations of Rost \cite{Ro-1} and Chacon
\cite{Chh}. He proved the existence of functions of time
(representing a reversed barrier) when $\mu$ has a continuous
distribution function which is flat around zero. He also showed that
these functions are unique under a Tychonov boundedness condition.

In this paper we develop an entirely different approach to the
embedding problem and prove the existence and uniqueness of
functions of time for general target laws $\mu$ with no extra
conditions imposed. The derivation of $\tau$ is constructive and the
construction itself is purely probabilistic and intuitive. The
method of proof relies upon weak convergence arguments for functions
of time arising from Helly's selection theorem and makes use of the
L\'evy metric which appears to be novel in the context of embedding
theorems. This enables us to avoid time-reversal arguments (present
in previous approaches) and relate the existence arguments directly
to the regularity of the sample path with respect to functions of
time. The fact that the construction applies to all target laws
$\mu$ with no integrability/regularity assumptions makes the
resulting embedding rather canonical and remarkable in the class of
known embeddings. Moreover, we show that the resulting stopping time
$\tau$ is minimal in the sense of Monroe \cite{Mo} so that the
stopped process $B^{\tau} = (B_{t \wedge\tau})_{t \ge0}$ satisfies
natural uniform integrability conditions which fail to hold for
trivial embeddings of any law (see, e.g., \cite{RY},
Exercise~5.7, page 276). We also show that the
resulting stopping time $\tau$ has
the smallest truncated expectation among all stopping times that
embed $\mu$ into $B$. The same result was derived by Chacon
\cite{Chh} for stopping times arising from the filling scheme when
their means are finite. A~converse result for stopping times arising
from barriers was first derived by Rost~\cite{Ro-2}. The main
results extend from standard Brownian motion to all recurrent
diffusion processes on the real line. Extending these results to
more general Markov processes satisfying specified conditions leads
to a research agenda which we leave open for future developments.

When the process is standard Brownian motion, then it is possible to
check that the sufficient conditions derived by Chacon (\cite{Chh}, page
47), are satisfied so that the filling scheme stopping time used
by Rost \cite{Ro-1} coincides with the first entry time to a
reversed barrier. If $\mu$ has a continuous distribution function
which is flat around zero, then the uniqueness result of McConnell
(\cite{Mc}, pages 684--690), implies that this reversed barrier is
uniquely determined under a Tychonov boundedness condition. When any
of these conditions fails, however, then it becomes unclear whether a
reversed barrier is uniquely determined by the filling scheme
because in principle there could be many reversed barriers yielding
the same law. One consequence of the present paper is that the
latter ambiguity gets removed since we show that the filling scheme
does indeed determine a reversed barrier uniquely for general target
laws $\mu$ with no extra conditions imposed. Despite this
contribution to the theory of filling schemes (see \cite{Chh} and the
references therein), it needs to be noted that the novel methodology
of the present paper avoids the filling scheme completely and
focuses on constructing the reversed barrier by functions of time
directly.

%s2 #&#
\section{Existence}\label{sec2}
%%%%%%%%%%%%%%%%%%%

In this section we state and prove the main existence result (see
also Corollary~\ref{cor8} below).

%th1 #&#
\begin{teo}[(Existence)]\label{teo1}
Let $B=(B_t)_{t \ge0}$ be a
standard Brownian motion defined on a probability space $(\Omega,
{\mathcal F}, \mathsf P)$ with $B_0=0$, and let $\mu$ be a probability
measure on
$(\mathbb{R},{\mathcal B}(\mathbb{R}))$ such that $\mu(\{0\})=0$.

\begin{longlist}[(III)]
\item[(I)] If $\operatorname{supp}(\mu) \subseteq\mathbb{R}_+$,
then there exists a
left-continuous increasing function $b\dvtx (0,\infty) \rightarrow
\mathbb{R}$
such that $B_{\tau_b} \sim\mu$ where $\tau_b = \inf\{ t>0
\vert B_t \ge b(t) \}$.
\end{longlist}

\begin{longlist}[(III)]
\item[(II)] If $\operatorname{supp}(\mu) \subseteq\mathbb{R}_-$,
then there exists a
left-continuous decreasing function $c\dvtx (0,\infty) \rightarrow
\mathbb{R}$
such that $B_{\tau_c} \sim\mu$ where $\tau_c = \inf\{ t>0
\vert B_t \le c(t) \}$.
\end{longlist}

\begin{longlist}[(III)]
\item[(III)] If $\operatorname{supp}(\mu) \cap\mathbb{R}_+ \ne
\varnothing$ and $\operatorname{supp}(\mu)
\cap\mathbb{R}_- \ne\varnothing$, then there exist a left-continuous
increasing function $b\dvtx (0,\infty) \rightarrow\mathbb{R}\cup\{
+\infty
\}$ and a left-continuous decreasing function $c\dvtx (0,\infty)
\rightarrow\mathbb{R}\cup\{ -\infty\}$ such that $B_{\tau_{b,c}}
\sim
\mu$ where $\tau_{b,c} = \inf\{ t>0 \vert B_t \ge b(t)
\mbox{ or } B_t \le c(t) \}$ (see Figure~\ref{fig1} below).
\end{longlist}
\end{teo}

\begin{pf} We will first derive (I)${}+{}$(II) since (III) will then
follow by combining and further extending the construction and
arguments of (I)${}+{}$(II). This will enable us to focus more clearly on
the subtle technical issues in relation to (a) the competing
character of the two boundaries in (III) and (b) the fact that one
of them can jump to infinity. Neither of these technical issues is
present in (I)${}+{}$(II) so that the key building block of the
construction is best understood by considering this case first.

%
%f1 #&#
\begin{figure}[t]

\includegraphics{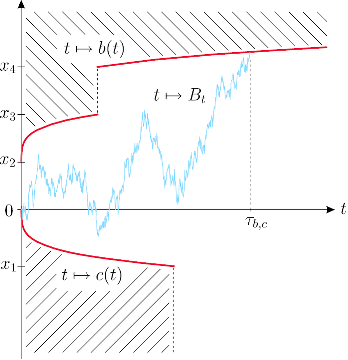}

\caption{An illustration of the reversed-barrier stopping
time $\tau_{b,c}$ from Theorem~\protect\ref{teo1} that embeds $\mu$ into $B$ when
$\operatorname{supp}(\mu) = [x_1,0] \cup[x_2,x_3] \cup[x_4,\infty)$.}\label{fig1}
\end{figure}

\begin{longlist}
\item[(I)${}+{}$(II) \textit{One-sided support}:] Clearly it is enough to prove
(I) since (II) then follows by symmetry. Let us therefore assume
that $\operatorname{supp}(\mu) \subseteq\mathbb{R}_+$\break throughout.
\end{longlist}

\begin{longlist}[1.]
\item[1. \textit{Bounded support}.]
Assume first that $\operatorname{supp}(\mu ) \subseteq
[0,\beta]$ for some $\beta< \infty$. Without loss of generality we
can assume that $\beta$ belongs to $\operatorname{supp}(\mu)$. Let
$0 = x_0^n <
x_1^n < \cdots< x_{m_n}^n = \beta$ be a partition of $[0,\beta]$
such that $\max_{ 1 \le k \le m_n} (x_k^n - x_{k-1}^n) \rightarrow
0$ as $n \rightarrow\infty$ (e.g., we could take a dyadic
partition defined by $x_k^n = \frac{k}{2^n} \beta$ for $k =
0,1, \ldots, 2^n$, but other choices are also possible and will lead
to the same result). Let $X$ be a random variable (defined on some
probability space) having the law equal to $\mu$, and set
%
%
%e2.1 #&#
\begin{equation}
\label{21} X_n = \sum_{k=1}^{m_n}
x_k^n I\bigl(x_{k-1}^n < X \le
x_k^n\bigr)
\end{equation}
for $n \ge1$. Then $X_n \rightarrow X$ almost surely, and hence $X_n
\rightarrow X$ in law as $n \rightarrow\infty$. Denoting the law of
$X_n$ by $\mu_n$, this means that $\mu_n \rightarrow\mu$ weakly as
$n \rightarrow\infty$. We will now construct a left-continuous
increasing function $b_n\dvtx (0,\infty) \rightarrow\mathbb{R}$ taking values
in $\{ x_1^n, x_2^n, \ldots, x_{m_n}^n \}$ such that $\tau_{b_n} =
\inf\{ t>0 \vert B_t \ge b_n(t) \}$ satisfies
\mbox{$B_{\tau_{b_n}} \sim\mu_n$} for~$n \ge1$.
\end{longlist}

\begin{longlist}[1.1.]
\item[1.1. \textit{Construction}: \textit{Discrete case}.] For this, set $p_k^n =
\mathsf P(x_{k-1}^n < X \le x_k^n)$ for $k=1,2, \ldots, m_n$ with
$n \ge1$ given and fixed, and let $k_1$ denote the smallest $k$ in
$\{ 1,2, \ldots, m_n \}$ such that $p_k^n>0$. Consider the
sequential movement of two sample paths $t \mapsto B_t$ and $t
\mapsto x_{k_1}^n$ as $t$ goes from $0$ onwards. From the recurrence
of $B$ it is clear that there exists a unique $t_1^n>0$ such that
the probability of $B$ hitting $x_{k_1}^n$ before $t_1^n$ equals
$p_{k_1}^n$. Stop the movement of $t \mapsto x_{k_1}^n$ at $t_1^n$,
and replace it with $t \mapsto x_{k_2}^n$ afterwards\vspace*{1pt} where $k_2$ is
the smallest $k$ in $\{ k_1 + 1, k_1 + 2, \ldots, m_n \}$
such that $p_k^n>0$. Set $b_n(t) = x_{k_1}^n$ for $t \in(0,t_1^n]$,
and on the event that $B$ did not hit $b_n$ on $(0,t_1^n]$, consider
the movement of $t \mapsto B_t$ and $t \mapsto x_{k_2}^n$ as $t$
goes from $t_1^n$ onwards. From the recurrence of $B$ it is clear
that there exists a unique $t_2^n > t_1^n$ such that the probability
of $B$ hitting $x_{k_2}^n$ before $t_2^n$ equals $p_{k_2}^n$.
Proceed as before, and set $b_n(t) = x_{k_2}^n$ for $t \in
(t_1^n,t_2^n]$. Continuing this construction by induction until
$t_i^n = \infty$ for some $i \le m_n$ (which clearly has to happen)
we obtain $b_n$ as stated above. Note that $b_n(t) = x_{k_1}^n$ for
$t \in(0,t_1^n]$ with $x_{k_1}^n \rightarrow\alpha=: \min
\operatorname{supp}(\mu)$ as $n \rightarrow\infty$ and $b_n(t) =
x_{m_n}^n$ for
$t \in(t_{i-1}^n,\infty)$ since $x_{m_n}^n = \beta= \max
\operatorname{supp}(\mu)$ by assumption.
\end{longlist}

\begin{longlist}[1.2.]
\item[1.2. \textit{Construction}: \textit{Passage to limit}.] In this way we have
obtained a sequence of left-continuous increasing functions $b_n\dvtx
(0,\infty) \rightarrow[\alpha,\beta]$ satisfying $b_n(0+)
\rightarrow\alpha$ as $n \rightarrow\infty$ and $b_n(+\infty) =
\beta$ for $n \ge1$. We can formally extend each $b_n$ to
$(-\infty,0]$ by setting $b_n(t) = b_n(0+)$ for $t \in(-1,0]$ and
$b_n(t) = 0$ for $t \in(-\infty,-1]$ (other definitions are also
possible). Then $\{ b_n \vert n \ge1 \}$ is a sequence of
left-continuous increasing functions from $\mathbb{R}$ into $\mathbb
{R}$ such that
$b_n(-\infty)=0$ and $b_n(+\infty) = \beta$ for all $n \ge1$. By
\emph{Helly's selection theorem} (see, e.g., \cite{Bi}, pages 336--337)
we therefore know that there exists a subsequence $\{ b_{n_k}
\vert k \ge1 \}$ and a left-continuous increasing function $b\dvtx
\mathbb{R}\rightarrow\mathbb{R}$ such that $b_{n_k} \rightarrow b$
weakly as $k
\rightarrow\infty$ in the sense that $b_{n_k}(t) \rightarrow b(t)$
as $k \rightarrow\infty$ for every $t \in\mathbb{R}$ at which $b$ is
continuous. (Note that since $b_n(t) = b_n(0+) \rightarrow\alpha$
as $n \rightarrow\infty$ for every $t \in(-1,0]$ it follows that
$b(0)=\alpha$ by the increase and left-continuity of $b$.)
Restricting $b$ to $(0,\infty)$ and considering the stopping time
%
%
%e2.2 #&#
\begin{equation}
\label{22} \tau_b = \inf\bigl\{ t>0 \vert B_t \ge
b(t) \bigr\},
\end{equation}
we claim that $B_{\tau_b} \sim\mu$. This can be seen as follows.
\end{longlist}

\begin{longlist}[1.3.]
\item[1.3. \textit{Tightness}.] We claim that the sequence of generalised
distribution functions $\{ b_n \vert n \ge1 \}$ is
\emph{tight} (in the sense the mass of the Lebesgue--Stieltjes
measure associated with $b_n$ cannot escape to infinity as $n
\rightarrow\infty$). Indeed, if $\varepsilon>0$ is given and fixed, then
$\delta_\varepsilon:= \mu((\beta- \varepsilon,\beta])
> 0$ since $\beta$ belongs to $\operatorname{supp}(\mu)$. Setting
$\tau_\beta=
\inf\{ t>0 \vert B_t \ge\beta\}$ we see that there
exists $t_\varepsilon>0$ large enough such that $\mathsf P(\tau
_\beta\le t_\varepsilon)
> 1 - \delta_\varepsilon$. Since $b_n \le\beta$ and hence $\tau
_{b_n} \le
\tau_\beta$ this implies that $\mathsf P(\tau_{b_n} \le
t_\varepsilon) > 1 -
\delta_\varepsilon$ for all $n \ge1$. From the construction of
$b_n$ the
latter inequality implies that $b_n(t_\varepsilon) > \beta-
\varepsilon$ for all
$n \ge1$. Recalling the extension of $b_n$ to $(-\infty,0]$
specified above where $b_n(-1)=0$, it therefore follows that
%
%
%e2.3 #&#
\begin{equation}
\label{23} b_n(t_\varepsilon) - b_n(-1) > \beta-
\varepsilon
\end{equation}
for all $n \ge1$. This shows that $\{ b_n \vert n \ge1 \}$
is tight as claimed. From (\ref{23}) we see that $b(+\infty) =
\beta$ and $b(-\infty)=0$ so that the Lebesgue--Stieltjes measure
associated with $b$ on $\mathbb{R}$ has a full mass equal to $\beta$ like
all other $b_n$ for $n \ge1$. Recalling that $b(0+)=\alpha$ we see
that the Lebesgue--Stieltjes measure associated with $b$ on
$(0,\infty)$ has a full mass equal to $\beta- \alpha$. For our
purposes we only need to consider the restriction of $b$ to
$(0,\infty)$.
\end{longlist}

\begin{longlist}[1.4.]
\item[1.4. \textit{L\'evy metric and convergence}.] If $b$ and $c$ are
left-continuous increasing functions from $\mathbb{R}$ into $\mathbb
{R}$ such that
$b(-\infty) = c(-\infty) = 0$ and $b(+\infty) =  c(+\infty) = \beta$,
then the \emph{L\'evy metric} is defined by
%
%
%e2.4 #&#
\begin{equation}
\label{24} \qquad d(b,c) = \inf\bigl\{ \varepsilon>0 \vert b(t - \varepsilon
) -\varepsilon\le c(t) \le b(t + \varepsilon) + \varepsilon\mbox{ for
all } t \in
\mathbb{R} \bigr\}.
\end{equation}
It is well known (see, e.g., \cite{Bi}, Exercise~14.5) that $c_n
\rightarrow b$ weakly if and only if $d(b,c_n) \rightarrow0$ as $n
\rightarrow\infty$. Defining functions
%
%
%e2.5 #&#
\begin{equation}
\label{25} b_\varepsilon(t):= b(t - \varepsilon) - \varepsilon\quad
\mbox{and}\quad b^\varepsilon(t):= b(t + \varepsilon) + \varepsilon
\end{equation}
for $t \in\mathbb{R}$, we claim that
%
%
%e2.6 #&#
%e2.7 #&#
\begin{eqnarray}
\label{26} \tau_{b_\varepsilon} &\uparrow&\tau_b\qquad \mathsf P\mbox{-a.s.},
\\
\label{27}  \tau_{b^\varepsilon} &\downarrow&\tau_b\qquad \mathsf P\mbox{-a.s.}
\end{eqnarray}
as $\varepsilon\downarrow0$, where in (\ref{26}) we also assume that
$b(0+)>0$.

\begin{pf*}{Proof of (\ref{26})}
Note first that $b_{\varepsilon'} \le
b_{\varepsilon''} \le b$ so that $\tau_{b_{\varepsilon'}} \le\tau
_{b_{\varepsilon''}}
\le\tau_b$ for $\varepsilon' \ge\varepsilon'' > 0$. It follows
therefore that
$\tau_{b_-}:= \lim_{ \varepsilon\downarrow0} \tau
_{b_\varepsilon} \le\tau_b$.
Moreover by definition of $\tau_{b_\varepsilon}$ we can find a sequence
$\delta_n \downarrow0$ as $n \rightarrow\infty$ such that $B_{
\tau_{b_\varepsilon}+\delta_n} \ge b_\varepsilon(\tau
_{b_\varepsilon} + \delta_n) = b(
\tau_{b_\varepsilon} - \varepsilon+ \delta_n) - \varepsilon$
for all $n \ge1$ with
$\varepsilon>0$. Letting $n \rightarrow\infty$ it follows that
$B_{\tau_{b_\varepsilon}} \ge b((\tau_{b_\varepsilon} -
\varepsilon)+) - \varepsilon\ge
b(\tau_{b_\varepsilon} - \varepsilon) - \varepsilon\ge b(\tau
_{b_\varepsilon} - \varepsilon_0) -
\varepsilon$ for all $\varepsilon\in(0, \varepsilon_0)$ with
$\varepsilon_0>0$ given and
fixed. Since $b$ is left-continuous and increasing, it follows that
$b$ is lower semicontinuous and hence by letting $\varepsilon
\downarrow0$
in the previous identity, we find that $B_{\tau_{b_-}} \ge
\liminf_{\varepsilon\downarrow0} (b(\tau_{b_\varepsilon} -
\varepsilon_0) - \varepsilon)
\ge b(\liminf_{\varepsilon\downarrow0} \tau_{b_\varepsilon} -
\varepsilon_0) =
b(\tau_{b_-} - \varepsilon_0)$ for all $\varepsilon_0>0$. Letting
$\varepsilon_0
\downarrow0$ and using that $b$ is left-continuous, we get
$B_{\tau_{b_-}} \ge b(\tau_{b_-})$. This implies that $\tau_b \le
\tau_{b_-}$, and hence $\tau_{b_-} = \tau_b$ as claimed in
(\ref{26}) above.
\end{pf*}

\begin{pf*}{Proof of (\ref{27})}
Note first that $b \le b^{\varepsilon
'} \le
b^{\varepsilon''}$ so that $\tau_b \le\tau_{b^{\varepsilon'}}
\le
\tau_{b^{\varepsilon''}}$ for $\varepsilon'' \ge\varepsilon' >
0$. It follows therefore
that $\tau_b \le\tau_{b^+}:= \lim_{ \varepsilon\downarrow0}
\tau_{b^\varepsilon}$. Moreover setting
%
%
%e2.8 #&#
\begin{equation}
\label{28} \sigma_b = \inf\bigl\{ t>0 \vert B_t >
b(t) \bigr\},
\end{equation}
we claim that
%
%
%e2.9 #&#
\begin{equation}
\label{29} \tau_b = \sigma_b\qquad\mathsf P\mbox{-a.s.}
\end{equation}
so that outside a $\mathsf P$-null set we have $B_{t_n} > b(t_n)$ for
some $t_n \downarrow\tau_b$ with $t_n > \tau_b$. Since $b$ is
increasing, each $t_n$ can be chosen as a continuity point of $b$,
and therefore there exists $\varepsilon_n>0$ small enough such that
$B_{t_n} > b^{\varepsilon_n}(t_n) = b(t_n + \varepsilon_n) +
\varepsilon_n > b(t_n)$
for all $n \ge1$. This shows that $\tau_{b^+} \le t_n$ outside the
$\mathsf P$-null set for all $n \ge1$. Letting $n \rightarrow\infty$
we get $\tau_{b^+} \le\tau_b$ $\mathsf P$-a.s. and hence $\tau
_{b^+} =
\tau_b$ $\mathsf P$-a.s. as claimed in (\ref{27}) above.
\end{pf*}

\begin{pf*}{Proof of (\ref{29})} Let us first introduce
%
%
%e2.10 #&#
\begin{equation}
\label{210} \tau_{b+\varepsilon} = \inf\bigl\{ t>0 \vert B_t \ge
b(t) + \varepsilon\bigr\},
\end{equation}
and note that $\tau_{b+}:= \lim_{ \varepsilon\downarrow0} \tau
_{b+\varepsilon}
= \sigma_b$ as is easily seen from definitions (\ref{28}) and~(\ref{210}). Next introduce the truncated versions of (\ref{22})
and (\ref{210}) by setting
%
%
%e2.11 #&#
%e2.12 #&#
\begin{eqnarray}
\label{211} \tau_b^\delta&=& \inf\bigl\{ t > \delta\vert
B_t \ge b(t) \bigr\},
\\
\label{212} \tau_{b+\varepsilon}^\delta&=& \inf\bigl\{ t > \delta
\vert B_t \ge b(t) + \varepsilon\bigr\}
\end{eqnarray}
with $\delta>0$ given and fixed. Note that $\tau_b^\delta\le
\tau_{b+\varepsilon'}^\delta\le\tau_{b+\varepsilon''}^\delta$
for $\varepsilon'' \ge
\varepsilon' > 0$. It follows therefore that $\tau_b^\delta\le
\tau_{b+}^\delta:= \lim_{ \varepsilon\downarrow0}
\tau_{b+\varepsilon}^\delta$. To prove that
%
%
%e2.13 #&#
\begin{equation}
\label{213} \tau_b^\delta= \tau_{b+}^\delta
\qquad\mathsf P\mbox{-a.s.}
\end{equation}
it is enough to establish that
%
%
%e2.14 #&#
\begin{equation}
\label{214} \mathsf P\bigl(\tau_{b+}^\delta> t\bigr) \le
\mathsf P\bigl( \tau_b^\delta> t\bigr)
\end{equation}
for all $t>0$. Indeed, in this case we have $\mathsf E (\tau
_{b+}^\delta
\wedge N) = \int_0^N \mathsf P(\tau_{b+}^\delta> t) \,dt \le\int_0^N
\mathsf P(\tau_b^\delta> t) \,dt = \mathsf E (\tau_b^\delta\wedge N)$
so that
$\tau_{b+}^\delta\wedge N = \tau_b^\delta\wedge N$ $\mathsf P$-a.s.
for all $N \ge1$. Letting $N \rightarrow\infty$ we obtain
(\ref{213}) as claimed. Assuming that (\ref{213}) is established,
note that
%
%
%e2.15 #&#
\begin{eqnarray}\label{215}
\sigma_b &=& \tau_{b+} =  \lim
_{\varepsilon\downarrow0} \tau_{b+\varepsilon} = \lim_{\varepsilon
\downarrow0} \lim
_{\delta\downarrow0} \tau_{b+
\varepsilon}^\delta= \lim
_{\delta\downarrow0} \lim_{\varepsilon
\downarrow
0} \tau_{b+\varepsilon}^\delta= \lim_{\delta\downarrow0} \tau_{b+} ^\delta
\nonumber\\[-8pt]\\[-8pt]
&=& \lim _{\delta\downarrow0} \tau_b^\delta= \tau_b
\qquad\mathsf P\mbox{-a.s.},\nonumber
\end{eqnarray}
where we use that $\varepsilon\mapsto\tau_{b+\varepsilon}^\delta
$ and $\delta
\mapsto\tau_{b+\varepsilon}^\delta$ are decreasing as
$\varepsilon\downarrow0$
and $\delta\downarrow0$ so that the two limits commute. Hence we
see that the proof of (\ref{29}) is reduced to establishing
(\ref{214}).
\end{pf*}

\begin{pf*}{Proof of (\ref{214})}
Note by Girsanov's theorem that
%
%
%e2.16 #&#
\begin{eqnarray}\label{216}
\mathsf P\bigl(\tau_{b+}^\delta> t\bigr)
&=& \mathsf P \Bigl( \lim_{\varepsilon\downarrow0} \tau_{b+\varepsilon
}^\delta>
t \Bigr) \le\lim_{\varepsilon
\downarrow0} \mathsf P\bigl( \tau_{b+\varepsilon}^\delta
> t \bigr)\nonumber
\\
& =& \lim_{\varepsilon
\downarrow0} \mathsf P\bigl(B_s < b(s) +
\varepsilon\mbox{ for all } s \in(\delta, t]\bigr)\nonumber
\\
&=& \lim_{\varepsilon\downarrow0} \mathsf P \biggl( B_s - \int
_0^s \frac{\varepsilon}{\delta} I(0 \le r \le\delta) \,dr
< b(s) \mbox{ for all } s \in(\delta, t] \biggr)
\\
&=& \lim_{\varepsilon\downarrow0} \mathsf E \biggl[ \frac
{{\mathcal E}_T^{H^\varepsilon}}{{\mathcal E}_T^{H^\varepsilon}} I \biggl(
B_s - \int_0^s H_r^\varepsilon \,dr < b(s) \mbox{ for all } s \in(\delta, t] \biggr) \biggr]\nonumber
\\
&=& \lim_{\varepsilon\downarrow0} \widetilde{\mathsf E} \biggl[
\frac{1} {{\mathcal E}_T^{H^\varepsilon}} I \bigl( \widetilde B_s < b(s)
\mbox{ for all } s \in
(\delta, t] \bigr) \biggr],\nonumber
\end{eqnarray}
where\vspace*{1pt} $H_r^\varepsilon= \frac{\varepsilon}{\delta} I(0 \le r
\le
\delta)$ and ${\mathcal E}_T^{H^\varepsilon} = \exp(\int_0^T
H_r^\varepsilon \,dB_r -
\frac{1}{2} \int_0^T (H_r^\varepsilon)^2 \,dr)$ so that~$d \widetilde
{\mathsf P}
= {\mathcal E}_T^{H^\varepsilon} \,d \mathsf P$ and $1/{\mathcal
E}_T^{H^\varepsilon} =
\exp(- \int_0^T H_r^\varepsilon \,dB_r + \frac{1}{2} \int_0^T
(H_r^\varepsilon)^2 \,dr) =\break \exp(- \int_0^T H_r^\varepsilon d \widetilde
B_r -
\frac{1}{2} \int_0^T (H_r^\varepsilon)^2 \,dr) =
\exp(-\frac{\varepsilon}{\delta} \widetilde B_\delta - \frac{1}{2}
\frac{\varepsilon^2}{\delta})$ with\vspace*{1pt} $\widetilde B_s = B_s - \int_0^s
H_r^\varepsilon \,dr$ being a standard Brownian motion under $\widetilde
{\mathsf P}$
for $s \in[0,T]$. From (\ref{216}) it therefore follows that
%
%
%e2.17 #&#
\begin{eqnarray}
\label{217} \qquad \mathsf P\bigl(\tau_{b+}^\delta> t\bigr) &\le&
\lim_{\varepsilon\downarrow
0} \mathsf E \biggl[ \exp\biggl( - \frac{\varepsilon}{\delta}
B_\delta- \frac{1}{2} \frac{\varepsilon^2}{\delta} \biggr) I \bigl(
B_s < b(s) \mbox{ for all } s \in(\delta, t] \bigr) \biggr]
\nonumber\\[-8pt]\\[-8pt]
& =& \mathsf P\bigl( B_s < b(s) \mbox{ for all } s \in(
\delta, t] \bigr) = \mathsf P\bigl(\tau_b^\delta> t\bigr)\nonumber
\end{eqnarray}
using the dominated convergence theorem since $\mathsf E e^{c \vert
B_\delta\vert} < \infty$ for $c>0$. This completes the verification
of (\ref{214}), and thus (\ref{27}) holds as well. [For a different
proof of (\ref{214}) in a more general setting, see the proof of
Corollary~\ref{cor8} below.]
\end{pf*}
\end{longlist}

\begin{longlist}[1.5.]
\item[1.5. \textit{Verification}.] To prove that $\tau_b$ from (\ref{22})
satisfies $B_{\tau_b} \sim\mu$, consider first the case when
$b(0+)>0$. Recall that $b_{n_k} \rightarrow b$ weakly and therefore
$d(b,b_{n_k}) \rightarrow0$ as $k \rightarrow\infty$ where $d$ is
the L\'evy metric defined in (\ref{24}). To simplify the notation
in the sequel, let us set $b_k:= b_{n_k}$ for $k \ge1$. This yields
the existence of $\varepsilon_k \downarrow0$ as $k \rightarrow
\infty$
such that $b_{\varepsilon_k}(t) \le b_k(t) \le b^{\varepsilon
_k}(t)$ for all $t>0$
and $k \ge1$ [recall that $b_{\varepsilon_k}$ and $b^{\varepsilon
_k}$ are defined
by (\ref{25}) above]. It follows therefore that $\tau
_{b_{\varepsilon_k}}
\le\tau_{b_k} \le\tau_{b^{\varepsilon_k}}$ for all $k \ge1$.
Letting $k
\rightarrow\infty$ and using (\ref{26}) and (\ref{27}) above, we
obtain $\tau_b = \lim_{ k \rightarrow\infty} \tau_{b_{\varepsilon
_k}} \le
\liminf_{ k \rightarrow\infty} \tau_{b_k} \le\limsup_{ k
\rightarrow\infty} \tau_{b_k} \le\lim_{ k \rightarrow\infty}
\tau_{b^{\varepsilon_k}} = \tau_b$ $\mathsf P$-a.s. This shows that
$\tau_b =
\lim_{ k \rightarrow\infty} \tau_{b_k}$ $\mathsf P$-a.s. and hence
$B_{\tau_b} = \lim_{ k \rightarrow\infty} B_{\tau_{b_k}}$
$\mathsf P$-a.s. Recalling that $B_{\tau_{b_k}} \sim\mu_k$ for $k
\ge1$
and that $\mu_k \rightarrow\mu$ weakly as $k \rightarrow\infty$, we
see that $B_{\tau_b} \sim\mu$ as claimed.

Consider next the case when $b(0+)=0$. With $\delta>0$ given and
fixed set $b^\delta:= b \vee\delta$ and $b_n^\delta:= b_n \vee
\delta$ for $n \ge1$. Since $b_k \rightarrow b$ weakly we see that
$b_k^\delta\rightarrow b^\delta$ weakly, and hence by the first part
of the proof above [since $b^\delta(0+) = \delta> 0$] we know
that $\tau_{b_k^\delta} \rightarrow\tau_{b^\delta}$ $\mathsf P$-a.s.
so that $B_{\tau_{b_k^\delta}} \rightarrow B_{\tau_{b^\delta}}$
$\mathsf P$-a.s. as $k \rightarrow\infty$. Moreover, since
$\tau_{b_k^\delta} \rightarrow\tau_{b_k}$ and $\tau_{b^\delta}
\rightarrow\tau_{b}$ as $\delta\downarrow0$ we see that
%
%
%e2.18 #&#
\begin{equation}
\label{218} B_{\tau_{b_k^\delta}} \rightarrow B_{\tau_{b_k}} \quad\mbox{and}
\quad B_{\tau_{b^\delta}} \rightarrow B_{\tau_{b}}
\end{equation}
as $\delta\downarrow0$. From the fact that the first convergence
in $\mathsf P$-probability is uniform over all $k \ge1$ in the sense
that we have
%
%
%e2.19 #&#
\begin{equation}
\label{219} \sup_{k \ge1} \mathsf P ( B_{\tau_{b_k^\delta}} \ne
B_{\tau_{b_k}} ) \le\sup_{k \ge1} \mu_k\bigl((0,
\delta]\bigr) \le\mu\bigl((0,\delta]\bigr) \rightarrow0
\end{equation}
as $\delta\downarrow0$, it follows that the limits in
$\mathsf P$-probability commute so that
%
%
%e2.20 #&#
\begin{equation}
\label{220} B_{\tau_b} = \lim_{\delta\downarrow0} B_{\tau_b^\delta}
= \lim_{\delta\downarrow0} \lim_{k \rightarrow\infty} B_{\tau
_{b_k}^\delta} =
\lim_{k \rightarrow\infty} \lim_{
\delta\downarrow0} B_{\tau_{b_k}^\delta} =
\lim_{k \rightarrow
\infty} B_{\tau_{b_k}}.
\end{equation}
Recalling again that $B_{\tau_{b_k}} \sim\mu_k$ for $k \ge1$ and
that $\mu_k \rightarrow\mu$ weakly as $k \rightarrow\infty$, we see
that $B_{\tau_b} \sim\mu$ in this case as well. Note also that the
same arguments show [by dropping the symbol $B$ from the left-hand
side of (\ref{219}) above] that $\tau_b = \lim_{ k \rightarrow
\infty} \tau_{b_k}$ in $\mathsf P$-probability. This will be used in the
proof of (III) below.
\end{longlist}

\begin{longlist}[2.]
\item[2. \textit{Unbounded support}.] Consider now the case when $\sup
\operatorname{supp}(\mu) = +\infty$. Let $X$ be a random
variable (defined on
some probability space) having the law equal to $\mu$, and set $X_N =
X \wedge\beta_N$ for some $\beta_N \uparrow\infty$ as $N
\rightarrow\infty$ such that $\mu((\beta_N - \varepsilon,\beta
_N])>0$ for
all $\varepsilon>0$ and $N \ge1$. Let $N \ge1$ be given and fixed.
Denoting the law of $X_N$ by $\mu_N$ we see that $\operatorname
{supp}(\mu_N)
\subseteq[0,\beta_N]$ with $\beta_N \in \operatorname{supp}(\mu_N)$. Hence by the
previous part of the proof we know that there exists a
left-continuous increasing function $b_N\dvtx (0,\infty) \rightarrow
\mathbb{R}$ such that $B_{\tau_{b_N}} \sim\mu_N$. Recall that this $b_N$
is obtained as the weak limit of a subsequence of the sequence of
simple functions constructed by partitioning $(0,\beta_N)$.
Extending the same construction to partitioning
$[\beta_N,\beta_{N+1})$ while keeping the obtained subsequence of
functions with values in $(0,\beta_N)$, we again know by the previous
part of the proof that there exists a left-continuous increasing
function $b_{N+1}\dvtx (0,\infty) \rightarrow\mathbb{R}$ such that
$B_{\tau_{b_{N+1}}} \sim\mu_{N+1}$. This $b_{N+1}$ is obtained as
the weak limit of a further subsequence of the previous subsequence
of simple functions. Setting $t_N = \inf\{ t>0 \vert b_N(t)
= \beta_N \}$ it is therefore clear that $b_{N+1}(t) = b_N(t)$ for
all $t \in(0,t_N]$. Continuing this process by induction and
noticing that $t_N \uparrow t_\infty$ as $N \rightarrow\infty$, we
obtain a function $b\dvtx (0,t_\infty) \rightarrow\mathbb{R}$ such
that $b(t)
= b_N(t)$ for all $t \in(0,t_N]$ and $N \ge1$. Clearly $b$ is
left-continuous and increasing since each $b_N$ satisfies these
properties. Moreover we claim that $t_\infty$ must be equal to
$+\infty$. For this, note that $\mathsf P(B_{\tau_b} \le x) =
\mathsf P(B_{\tau_{b_N}} \le x)$ for $x < \beta_N$ and $N \ge1$.
Letting $N \rightarrow\infty$ and using that $B_{\tau_{b_N}} \sim
\mu_N$ converges weakly to $\mu$ since $X_N \rightarrow X$, we see
that $\mathsf P(B_{\tau_b} \le x) = \mathsf P(X \le x)$ for all $x>0$
at which
the distribution function of $X$ is continuous. Letting $x \uparrow
\infty$ over such continuity points we get $\mathsf P(B_{\tau_b} <
\infty) = 1$. Since clearly this is not possible if $t_\infty$ is
finite, we see that $t_\infty= +\infty$ as claimed. Noting that $b_N
= b \wedge\beta_N$ on $(0,\infty)$ for $N \ge1$ it follows that~$\tau_{b_N} = \inf\{ t>0 \vert B_t \ge b_N(t) \} = \inf
\{ t>0 \vert B_t \ge b(t) \wedge\beta_N \}$ from where we
see that $\tau_{b_N} \rightarrow\tau_b$ and thus
$B_{\tau_{b_N}} \rightarrow B_{\tau_b}$ as $N \rightarrow\infty$.
Since $X_N \rightarrow X$ and thus $\mu_N \rightarrow\mu$ weakly as
$N \rightarrow\infty$, it follows that $B_{\tau_b} \sim\mu$ as
claimed. This completes the proof of~(I).
\end{longlist}

\begin{longlist}[(III)]
\item[(III) \textit{Two-sided support}:] This will be proved by combining
and further extending the construction and arguments of (I) and
(II). Novel aspects in this process include the competing character
of the two boundaries and the fact that one of them can jump to
infinite value.
\end{longlist}

\begin{longlist}[3.]
\item[3. \textit{Bounded support}.] As in the one-sided case assume first
that $\operatorname{supp}(\mu) \subseteq[\gamma,\beta]$ for some
$\gamma<0<\beta$. Without loss of generality we can assume that
$\beta$ and $\gamma$ belong to $\operatorname{supp}(\mu)$. Let $0 =
x_0^n < x_1^n
< \cdots< x_{m_n}^n = \beta$ be a partition of $[0,\beta]$ such
that $\max_{ 1 \le k \le m_n} (x_k^n - x_{k-1}^n) \rightarrow0$
as $n \rightarrow\infty$, and let $0 = y_0^n > y_1^n
> \cdots> y_{l_n}^n = \gamma$ be a partition of $[\gamma,0]$ such
that\vspace*{1pt} $\max_{ 1 \le j \le l_n} (y_{j-1}^n - y_j^n) \rightarrow0$
as $n \rightarrow\infty$. Let $X$ be a random variable (defined on
some probability space) having the law equal to $\mu$, and set
%
%
%e2.21 #&#
\begin{equation}
\label{221} \qquad\quad X_n^+ = \sum_{k=1}^{m_n}
x_k^n I\bigl(x_{k-1}^n < X \le
x_k^n\bigr) \quad\mbox{and}\quad X_n^- =
\sum_{j=1}^{l_n} y_j^n
I\bigl(y_j^n \le X < y_{j-1}^n
\bigr)
\end{equation}
for $n \ge1$. Then $X_n^+ + X_n^- \rightarrow X$ almost surely and
hence $X_n^+ + X_n^- \rightarrow X$ in law as $n \rightarrow
\infty$. Denoting the law of $X_n^+ + X_n^-$ by $\mu_n$ and
recalling that $X$ has the law $\mu$, this means that $\mu_n
\rightarrow\mu$ weakly as $n \rightarrow\infty$. We will now
construct a left-continuous increasing function $b_n\dvtx (0,\infty)
\rightarrow\mathbb{R}$ taking values in $\{ x_1^n, x_2^n, \ldots,
x_{m_n}^n, +\infty\}$ and a left-continuous decreasing function
$c_n\dvtx (0,\infty) \rightarrow\mathbb{R}$ taking values in $\{
y_1^n, y_2^n,
\ldots, y_{l_n}^n, -\infty\}$ with $b_n(t) < +\infty$ or $c_n(t) >
-\infty$ for all $t \in(0,\infty)$ such that $\tau_{b_n,c_n} =
\inf\{ t>0 \vert B_t \ge b_n(t) \mbox{ or } B_t \le
c_n(t) \}$ satisfies $B_{\tau_{b_n,c_n}} \sim\mu_n$ for $n \ge
1$.
\end{longlist}

\begin{longlist}[3.1.]
\item[3.1. \textit{Construction}: \textit{Discrete case}.] For this, set $p_k^n =
\mathsf P(x_{k-1}^n < X \le x_k^n)$ for $k=1,2, \ldots, m_n$ and
$q_j^n = \mathsf P(y_j^n \le X < y_{j-1}^n)$ for $j=1,2, \ldots,
l_n$ with $n \ge1$ given and fixed. Let $k_1$ denote the smallest
$k$ in $\{ 1,2, \ldots, m_n \}$ such that $p_k^n > 0$, and let $j_1$
denote the smallest $j$ in $\{ 1,2, \ldots, l_n \}$ such that $q_j^n
> 0$. Consider the sequential movement of three sample paths $t
\mapsto B_t$, $t \mapsto x_{k_1}^n$ and $t \mapsto y_{j_1}^n$ as $t$
goes from $0$ onwards. From the recurrence of $B$ it is clear that
there exists a unique $t_1^n > 0$ such that the probability of $B$
hitting $x_{k_1}^n$ before $y_{j_1}^n$ on $(0,t_1^n]$ equals
$p_{k_1}^n$, or the probability of $B$ hitting $y_{j_1}^n$ before
$x_{k_1}^n$ on $(0,t_1^n]$ equals $q_{j_1}^n$, whichever happens
first (including simultaneous happening). In the first case stop the
movement of $t \mapsto x_{k_1}^n$ at $t_1^n$ and replace it with $t
\mapsto x_{k_2}^n$ afterwards where $k_2$ is the smallest $k$ in $\{
k_1 + 1,k_1 + 2, \ldots, m_n \}$ such that $p_k^n > 0$ (if there
is no such $k$ then make no further replacement). In the second case
stop the movement of $t \mapsto y_{j_1}^n$ at $t_1^n$, and replace it
with $t \mapsto y_{j_2}^n$ afterwards where $j_2$ is the smallest
$j$ in $\{ j_1 + 1,j_1 + 2, \ldots, l_n \}$ such that $q_j^n > 0$
(if there is no such $j$ then make no further replacement). In the
third case, when the first and second case happen simultaneously,
stop the movement of both $t \mapsto x_{k_1}^n$ and $t \mapsto
y_{j_1}^n$ at $t_1^n$, and replace them with $t \mapsto x_{k_2}^n$
and $t \mapsto y_{j_2}^n$, respectively (if there is no $k_2$ or
$j_2$, then make no replacement, resp.). In all three cases set
$b_n(t) = x_{k_1}^n$ and $c_n(t) = y_{j_1}^n$ for $t \in(0,t_1^n]$.
On the event that $B$ did not hit $b_n$ or $c_n$ on $(0,t_1^n]$, in
the first case consider the movement\vspace*{1pt} of $t \mapsto B_t$, $t \mapsto
x_{k_2}^n$ and $t \mapsto y_{j_1}^n$, in the second case consider
the movement of $t \mapsto B_t$, $t \mapsto x_{k_1}^n$ and $t
\mapsto y_{j_2}^n$, and in the third case\vspace*{1pt} consider the movement of
$t \mapsto B_t$, $t \mapsto x_{k_2}^n$, and $t \mapsto y_{j_2}^n$ as
$t$ goes from $t_1^n$ onwards. If there is no $k_2$ or $j_2$ we can
formally set $x_{k_2}^n = +\infty$ or $y_{j_2}^n = -\infty$,
respectively (note, however, that either $k_2$ or $j_2$ will always be
finite). Continuing this construction by induction until $t_i^n =
\infty$ for some $i \le m_n \vee l_n$ (which clearly has to happen)
we obtain $b_n$ and $c_n$ as stated above.
\end{longlist}

\begin{longlist}[3.2.]
\item[3.2. \textit{Construction}: \textit{Passage to limit}.] For $n \ge1$ given
and fixed note that $b_n$ takes value $\beta$ on some interval, and
$c_n$ takes value $\gamma$ on some interval since both $\beta$ and
$\gamma$ belong to $\operatorname{supp}(\mu)$. The main technical
difficulty is
that either $b_n$ can take value $+\infty$ or $c_n$ can take value
$-\infty$ from some time $t_\zeta$ onwards as well (in which case
the corresponding interval is bounded). In effect this means that
the corresponding function is not defined on $(t_\zeta,\infty)$ with
values in $\mathbb{R}$. To overcome this difficulty we will set $\bar b_n(t)
= \beta$ and $\bar c_n(t) = \gamma$ for $t>t_\zeta$. Setting further
$\bar b_n = b_n$ and $\bar c_n = c_n$ on $(0,t_\zeta]$ we see that
$\bar b_n$ and $\bar c_n$ are generalised distribution functions on
$(0,\infty)$. Note that we always have either $\bar b_n = b_n$ or
$\bar c_n = c_n$ (and often both). Note also that $\bar b_n \ne b_n$
if and only if $b_n$ takes value $+\infty$ and $\bar c_n \ne c_n$ if
and only if $c_n$ takes value $-\infty$. Note finally that $\bar
b_n(+\infty)=\beta$ and $\bar c_n(+\infty)=\gamma$. Applying the
same arguments as in Part~1.2 above (upon extending $\bar b_n$ and
$\bar c_n$ to $\mathbb{R}$ first) we know that there exist subsequences
$\{ \bar b_{n_k} \vert k \ge1 \}$ and $\{ \bar c_{n_k}
\vert k \ge1 \}$ such that $\bar b_{n_k} \rightarrow\bar b$
and $\bar c_{n_k} \rightarrow\bar c$ weakly as $k \rightarrow
\infty$ for some increasing left-continuous function $\bar b$ and
some decreasing left-continuous function $\bar c$.
\end{longlist}

\begin{longlist}[3.3.]
\item[3.3. \textit{Tightness}.] We claim that the sequences of generalised
distribution functions $\{ \bar b_n \vert n \ge1 \}$ and
$\{ \bar c_n \vert n \ge1 \}$ are \emph{tight}. Indeed, if
$\varepsilon>0$ is given and fixed, then $\delta_\varepsilon':=
\mu((\beta-
\varepsilon,\beta])>0$ and $\delta_\varepsilon'':= \mu([\gamma, \gamma+
\varepsilon))>0$ since $\beta$ and $\gamma$ belong to
$\operatorname{supp}(\mu)$. Setting
$\delta_\varepsilon:= \delta_\varepsilon' \wedge\delta
_\varepsilon''$ and considering
$\tau_\beta= \inf\{ t>0 \vert B_t \ge\beta\}$ and
$\tau_\gamma= \inf\{ t>0 \vert B_t \le\gamma\}$, we see
that there exists $t_\varepsilon>0$ large enough such that $\mathsf
P(\tau_\beta
\vee\tau_\gamma\le t_\varepsilon) > 1 - \delta_\varepsilon$. Since
$\tau_{b_n,c_n} \le\tau_\beta\vee\tau_\gamma$, this implies that
$\mathsf P(\tau_{b_n,c_n} \le t_\varepsilon)
> 1 - \delta_\varepsilon$ for all $n \ge1$. From the construction of
$b_n$ and $c_n$ the latter inequality implies that $b_n(t_\varepsilon
) >
\beta- \varepsilon$ and $c_n(t_\varepsilon) < \gamma+
\varepsilon$ for all $n \ge1$
(note that in all these arguments we can indeed use unbarred
functions). The tightness claim then follows using the same
arguments as in Part~1.3 above.
\end{longlist}

\begin{longlist}[3.4.]
\item[3.4. \textit{Verification}.] Applying the same arguments as in Part~1.4 above we know from Part~1.5 above that setting $\bar b_k:= \bar
b_{n_k}$ for $k \ge1$, we have $\tau_{\bar b_k} \rightarrow
\tau_{\bar b}$ and $\tau_{\bar c_k} \rightarrow\tau_{\bar c}$ in
$\mathsf P$-probability as $k \rightarrow\infty$. Setting $t_k^b =
\sup\{ t>0 \vert b_k(t)=\beta\}$ and $t_k^c = \sup
\{ t>0 \vert c_k(t)=\gamma\}$ by the construction above, we
know that either $t_k^b = \infty$ or $t_k^c = \infty$ for all $k \ge
1$. If there exists $k_0 \ge1$ such that both $t_k^b = \infty$ and
$t_k^c = \infty$ for all $k \ge k_0$, then $b_k = \bar b_k$ and $c_k
= \bar c_k$ for all $k \ge k_0$ so that $\tau_{b_k,c_k} = \tau_{b_k}
\wedge\tau_{c_k} = \tau_{\bar b_k} \wedge\tau_{\bar c_k}
\rightarrow\tau_{\bar b} \wedge\tau_{\bar c} = \tau_{\bar b,\bar c} =
\tau_{b,c}$ in $\mathsf P$-probability as $k \rightarrow\infty$
where we set $b:= \bar b$ and $c:= \bar c$. This implies that
$B_{\tau_{b_k,c_k}} \rightarrow B_{\tau_{b,c}}$ in
$\mathsf P$-probability and thus in law as well while
$B_{\tau_{b_k,c_k}} \sim\mu_k$ with $\mu_k \rightarrow\mu$
weakly as $k \rightarrow\infty$ then shows that
$B_{\tau_b,\tau_c} \sim\mu$ as required. Suppose therefore that
there is no such $k_0 \ge1$. This means that we have infinitely
many $t_k^b < \infty$ or infinitely many $t_k^c < \infty$ for $k \ge
1$. Without loss of generality assume that the former holds. Then we
can pass to a further subsequence such that $t_{k_l}^b < \infty$ for
all $l \ge1$ and $t_{k_l}^b \rightarrow t_\infty^b \in(0,\infty]$
as $l \rightarrow\infty$. Set $b(t) = \bar b(t)$ for $t \in
(0,t_\infty^b]$ and $b(t) = \infty$ for $t \in(t_\infty^b,\infty)$.
Set also $c(t) = \bar c(t)$ for $t>0$, and note that $c_{k_l} = \bar
c_{k_l}$ for all $l \ge1$. To simplify the notation set further
$b_l:= b_{k_l}$, $\bar b_l:= \bar b_{k_l}$, $c_l:= c_{k_l}$, and
$\bar c_l:= \bar c_{k_l}$ for $l \ge1$. Then $\tau_{\bar b_l}
\rightarrow\tau_{\bar b}$ in $\mathsf P$-probability and hence
$\tau_{\bar b_l} I(\tau_{\bar b} < t_\infty^b) \rightarrow
\tau_{\bar b} I(\tau_{\bar b} < t_\infty^b)$ in
$\mathsf P$-probability as $l \rightarrow\infty$. Using definitions of
barred functions and the fact that $t_{k_l}^b \rightarrow
t_\infty^b$, one can easily verify that the previous relation implies
that $\tau_{b_l} I(\tau_{b} < t_\infty^b) \rightarrow\tau_b
I(\tau_b < t_\infty^b)$ in $\mathsf P$-probability as $l \rightarrow
\infty$. Since $\mathsf P(\tau_b < t_\infty^b)=1$ it follows that
$\tau_{b_l} \wedge\tau_{c_l} \rightarrow\tau_b \wedge\tau_c$ in
$\mathsf P$-probability as $l \rightarrow\infty$. This implies that
$B_{\tau_{b_l,c_l}} \rightarrow B_{\tau_{b,c}}$ in
$\mathsf P$-probability as $l \rightarrow\infty$ and hence
$B_{\tau_{b,c}} \sim\mu$ using the same argument as above. This
completes the proof in the case when $\operatorname{supp}(\mu)$ is bounded.
\end{longlist}

\begin{longlist}[4.]
\item[4. \textit{Half bounded support}.] Consider now the case when $\sup
\operatorname{supp}(\mu) = +\infty$ and $\inf\operatorname
{supp}(\mu) =: \gamma\in
(-\infty,0)$; see Figure~\ref{fig1} above. Let $X$ be a random variable
(defined on some probability space) having the law equal to $\mu$,
and set $X_N = X \wedge\beta_N$ for some $\beta_N \uparrow\infty$
as $N \rightarrow\infty$ such that $\mu((\beta_N -
\varepsilon,\beta_N])>0$ for all $\varepsilon>0$ and $N \ge1$.
Let $N \ge1$ be
given and fixed. Denoting the law of $X_N$ by $\mu_N$ we see that
$\operatorname{supp}(\mu_N) \subseteq[\gamma,\beta_N]$ with $\beta
_N$ and
$\gamma$ belonging to $\operatorname{supp}(\mu_N)$. Hence by Parts
3.1--3.4 above
we know that there exist a left-continuous increasing function
$b_N\dvtx (0,\infty) \rightarrow(0,\beta_N] \cup\{+\infty\}$ and a
left-continuous decreasing function $c_N\dvtx (0,\infty) \rightarrow
[\gamma,0) \cup\{-\infty\}$ such that $B_{\tau_{b_N,c_N}} \sim
\mu_N$.
\end{longlist}

\begin{longlist}[4.1.]
\item[4.1. \textit{Construction}.] Recall that these $b_N$ and $c_N$ are
obtained as the weak limits of subsequences of the sequences of
simple functions constructed by partitioning $(\gamma,0)$ and
$(0,\beta_N)$. Extending the same construction to partitioning
$(\gamma,0)$ and $[\beta_N,\beta_{N+1})$ while keeping the obtained
subsequence of functions with values strictly smaller than $\beta_N$,
we again know by Parts 3.1--3.4 above that there exist a
left-continuous increasing function $b_{N+1}\dvtx (0,\infty)
\rightarrow(0,\beta_{N+1}] \cup\{+\infty\}$ and a left-continuous
decreasing function $c_{N+1}\dvtx (0,\infty) \rightarrow[\gamma,0)
\cup\{-\infty\}$ such that $B_{\tau_{b_{N+1},c_{N+1}}} \sim
\mu_{N+1}$. These $b_{N+1}$ and $c_{N+1}$ are obtained as the weak
limits of further subsequences of the previous subsequences of
simple functions. Setting $t_N = \inf\{ t>0 \vert b_N(t) =
\beta_N \}$ it is therefore clear that $b_{N+1}(t) = b_N(t)$ and
$c_{N+1}(t) = c_N(t)$ for all $t \in(0,t_N]$. Continuing this
process by induction and noticing that $t_N \uparrow t_\infty$ as $N
\rightarrow\infty$, we obtain a left-continuous increasing function
$b\dvtx (0,t_\infty) \rightarrow\mathbb{R}$ and a left-continuous decreasing
$c\dvtx (0,t_\infty) \rightarrow\mathbb{R}\cup\{-\infty\}$ such that
$b(t) =
b_N(t)$ and $c(t) = c_N(t)$ for all $t \in(0,t_N]$ and $N \ge1$.
Note that $b$ is finite valued on $(0,t_\infty)$ with $b(t_\infty-)
= +\infty$.
\end{longlist}

\begin{longlist}[4.2.]
\item[4.2. \textit{Verification}.] To verify that $b$ and $c$ are the
required functions, consider first the case when $t_\infty= \infty$.
If $c$ is finite valued, then $\tau_{b,c} < \infty$ $\mathsf P$-a.s. and
hence $\tau_{b_N,c_N} \rightarrow\tau_{b,c}$ $\mathsf P$-a.s. as $N
\rightarrow\infty$. If $c$ is not finite valued, then $c=c_N$ and
hence $\mathsf P(B_{\tau_{b,c}} < \beta_N) = \mathsf P(B_{\tau
_{b_N,c_N}} <
\beta_N) = 1 - \mu([\beta_N,\infty))$ for all $N \ge N_0$ with some
$N_0 \ge1$. Letting $N \rightarrow\infty$ and using that
$\mu([\beta_N,\infty)) \rightarrow0$, we find that $\mathsf P(\tau
_{b,c} <
\infty)=1$ and hence $\tau_{b_N,c_N} \rightarrow\tau_{b,c}$
$\mathsf P$-a.s. Thus the latter relation always holds and hence
$B_{\tau_{b_N,c_N}} \rightarrow B_{\tau_{b,c}}$ $\mathsf P$-a.s. as $N
\rightarrow\infty$. Since $B_{\tau_{b_N,c_N}} \sim\mu_N$ and
$X_N \rightarrow X$ so that $\mu_N \rightarrow\mu$ weakly as $N
\rightarrow\infty$ it follows that $B_{\tau_{b,c}} \sim\mu$ as
required.

Consider next the case when $t_\infty< \infty$. To extend the
function $c$ to $[t_\infty,\infty)$ when $c(t_\infty-) > \gamma$
(note that when $c(t_\infty-) = \gamma$ then clearly $c$ must remain
equal to $\gamma$ on $[t_\infty,\infty)$ as well) set $t_N^c =
\sup\{ t>0 \vert c_N(t) = \gamma\}$ and define $\bar c_N(t) = c_N(t)$
for $t \in(0,t_N^c]$ and $\bar c_N(t) = \gamma$
for $t \in(t_N^c,\infty)$ whenever $t_N^c < \infty$ for $N \ge1$.
Applying the same arguments as in Parts 1.2 and 1.3 above, we know
that there exists a subsequence $\{ \bar c_{N_k} \vert k \ge
1 \}$ and a left-continuous function $\bar c$ such that $\bar c_{N_k}
\rightarrow\bar c$ weakly as $k \rightarrow\infty$.
Applying the same arguments as in Part~1.4 above we know from Part~1.5 above that setting $\bar c_k:= \bar c_{N_k}$ for $k \ge1$ we
have $\tau_{\bar c_k} \rightarrow\tau_{\bar c}$ in
\mbox{$\mathsf P$-}probability as $k \rightarrow\infty$. Moreover, we claim
that $t_N^c \rightarrow\infty$ as $N \rightarrow\infty$. For this,
suppose that $t_{N_l}^c \le T < \infty$ for $l \ge1$. Fix
$\varepsilon>
0$ small and set $c_\varepsilon(t) = c(t)$ for $t \in(0,t_\infty-
\varepsilon)$
and $c_\varepsilon(t) = c(t_\infty- \varepsilon)$ for $t \in
[t_\infty-
\varepsilon,T]$. Setting $b_l:= b_{N_l}$ and $c_l:= c_{N_l}$ we
then have
$\mu([\gamma,\beta_{N_l})) = \mathsf P(B_{\tau_{b_l,c_l}} \in
[\gamma,\beta_{N_l})) \le\mathsf P(\tau_{b,c_\varepsilon} \le T)$
for all $l \ge
1$. Letting $l \rightarrow\infty$ and using that
$\mu([\gamma,\beta_{N_l})) \rightarrow1$, we see that
$\mathsf P(\tau_{b,c_\varepsilon} \le T) = 1$ which clearly is
impossible since
$b$ is not defined beyond $t_\infty$. Thus $t_N^c \rightarrow
\infty$ as $N \rightarrow\infty$ and hence $t_{N_k}^c \rightarrow
\infty$ as $k \rightarrow\infty$. Setting $c:= \bar c$ and $c_k:=
c_{N_k}$ for $k \ge1$ and using the same arguments as in Part~3.4
above, we can therefore conclude that $\tau_{c_k} I(\tau_c <
\infty) \rightarrow\tau_c I(\tau_c < \infty)$ in
$\mathsf P$-probability as $k \rightarrow\infty$. Since $\mathsf
P(\tau_c <
\infty)=1$ this shows that $\tau_{c_k} \rightarrow\tau_c$ in
$\mathsf P$-probability as $k \rightarrow\infty$. Setting $b_k:=
b_{N_k}$ and noting that $\tau_{b_k} \rightarrow\tau_b$ on $\{
\tau_b < \infty\}$, we see that $\tau_{b_k,c_k} \rightarrow
\tau_{b,c}$ in $\mathsf P$-probability as $k \rightarrow\infty$ and
hence $B_{\tau_{b,c}} \sim\mu$ using the same argument as above.
The case when $\sup\operatorname{supp}(\mu) \in(0,+\infty)$ and
$\inf
\operatorname{supp}(\mu) = -\infty$ follows in exactly the same way
by symmetry.
\end{longlist}

\begin{longlist}[5.]
\item[5. \textit{Fully unbounded support}.] Consider finally the remaining case when
both $\sup\operatorname{supp}(\mu) = +\infty$ and $\inf
\operatorname{supp}(\mu) =
-\infty$. Let $X$ be a random variable (defined on some probability
space) having the law equal to $\mu$, and set $X_N = \gamma_N \vee X
\wedge\beta_N$ for some $\beta_N \uparrow\infty$ and $\gamma_N
\downarrow-\infty$ as $N \rightarrow\infty$ such that
$\mu((\beta_N - \varepsilon,\beta_N])>0$ and $\mu([\gamma_N,
\gamma_N +
\varepsilon))>0$ for all $\varepsilon>0$ and $N \ge1$. Let $N \ge
1$ be given and
fixed. Denoting the law of $X_N$ by $\mu_N$ we see that
$\operatorname{supp}(\mu_N) \subseteq[\gamma_N,\beta_N]$ with
$\beta_N$ and
$\gamma_N$ belonging to $\operatorname{supp}(\mu_N)$. Hence by Parts
3.1--3.4 above
we know that there exist a left-continuous increasing function
$b_N\dvtx (0,\infty) \rightarrow(0,\beta_N] \cup\{+\infty\}$ and a
left-continuous decreasing function $c_N\dvtx (0,\infty) \rightarrow
[\gamma_N,0) \cup\{-\infty\}$ such that $B_{\tau_{b_N,c_N}} \sim
\mu_N$.
\end{longlist}

\begin{longlist}[5.1.]
\item[5.1. \textit{Construction}.] Recall that these $b_N$ and $c_N$ are
obtained as the weak limits of subsequences of the sequences of
simple functions constructed by partitioning $(\gamma_N,0)$ and
$(0,\beta_N)$. Extending the same construction to partitioning
$(\gamma_{N+1},\gamma_N]$ and $[\beta_N,\beta_{N+1})$ while keeping
the obtained subsequence of functions with values strictly smaller
than $\beta_N$ and strictly larger than $\gamma_N$, we again know by
Parts 3.1--3.4 above that there exist a left-continuous increasing
function $b_{N+1}\dvtx (0,\infty) \rightarrow(0,\beta_{N+1}] \cup
\{+\infty\}$ and a left-continuous decreasing function $c_{N+1}\dvtx
(0,\infty) \rightarrow[\gamma_{N+1},0) \cup\{-\infty\}$ such that
$B_{\tau_{b_{N+1},c_{N+1}}} \sim\mu_{N+1}$. These $b_{N+1}$ and
$c_{N+1}$ are obtained as the weak limits of further subsequences of
the previous subsequences of simple functions. Setting $t_N^b =
\inf\{ t>0 \vert b_N(t) = \beta_N \}$ and $t_N^c = \inf
\{ t>0 \vert c_N(t) = \gamma_N \}$ it is therefore clear
that $b_{N+1}(t) = b_N(t)$ and $c_{N+1}(t) = c_N(t)$ for all $t \in
(0,t_N]$ where we set $t_N:= t_N^b \wedge t_N^c$ for $N \ge1$.
Continuing this process by induction and noticing that $t_N \uparrow
t_\infty$ as $N \rightarrow\infty$, we obtain a left-continuous
increasing function $b\dvtx (0,t_\infty) \rightarrow\mathbb{R}$ and a
left-continuous decreasing $c\dvtx (0,t_\infty) \rightarrow\mathbb{R}$ such
that $b(t) = b_N(t)$ and $c(t) = c_N(t)$ for all $t \in(0,t_N]$ and
$N \ge1$.
\end{longlist}

\begin{longlist}[5.2.]
\item[5.2. \textit{Verification}.] To verify that $b$ and $c$ are the
required functions, consider first the case when $t_\infty= \infty$.
Then since $b(t_N) \le\beta_N$ and $c(t_N) \ge\gamma_N$ for any $A
\in{\mathcal B}(\mathbb{R})$, we have $\mathsf P(B_{\tau_{b,c}} \in
A \cap
(c(t_N),b(t_N))) = \mathsf P(B_{\tau_{b_N,c_N}} \in A \cap
(c(t_N),b(t_N))) = \mu(A \cap(c(t_N),b(t_N)))$ for all $N \ge1$.
Letting $N \rightarrow\infty$ and using that $b(t_N) \uparrow
\infty$ and $c(t_N) \downarrow-\infty$, we see that
$\mathsf P(B_{\tau_{b,c}} \in A) = \mu(A)$, and this shows that
$B_{\tau_{b,c}} \sim\mu$ as required.

Consider next the case when $t_\infty< \infty$, and assume first
that either $\{ t_N^b \vert N \ge1 \}$ or $\{ t_N^c
\vert N \ge1 \}$ is not bounded (we will see below that this is
always true). Without loss of generality we can assume (by passing
to a subsequence if needed) that $t_N^c \rightarrow\infty$ so that
$t_N^b \uparrow t_\infty< \infty$ as $N \rightarrow\infty$. To
extend the function $c$ to $[t_\infty,\infty)$ we can now connect to
the final paragraph of Part~4 above. Choosing $M \ge1$ large enough
so that $\gamma_M < c(t_\infty-)$, we see that we are in the setting
of that paragraph with $\gamma= \gamma_M$, and hence there exists a
left-continuous decreasing function $c_M\dvtx (0,\infty) \rightarrow
[\gamma_M,0)$ such that $B_{\tau_{b,c_M}} \sim X \vee\gamma_M$.
Recall that this $c_M$ is obtained as the weak limit of a
subsequence of the sequence of functions embedding $B$ into
$[\gamma_M,\beta_N]$ for $N \ge1$, and note that $c_M$ coincides
with $c$ on $(0,t_\infty)$. Extending the same construction to
embedding $B$ into $[\gamma_{M+1},\beta_N]$ for $N \ge1$ while
keeping the subsequence of functions obtained previously, we again
know by the final paragraph of Part~4 above that there exists a
left-continuous decreasing function $c_{M+1}\dvtx (0,\infty)
\rightarrow[\gamma_{M+1},0)$ such that $B_{\tau_{b,c_{M+1}}} \sim
X \vee\gamma_{M+1}$. This $c_{M+1}$ is obtained as the weak limit
of a further subsequence of the previous sequence of functions.
Setting $t_M^c = \inf\{ t>0 \vert c_M(t) = \gamma_M \}$ it
is therefore clear that $c_{M+1}(t) = c_M(t)$ for $t \in(0,t_M^c)$.
Continuing this process by induction we obtain a left-continuous
decreasing function $c\dvtx (0,\infty) \rightarrow\mathbb{R}$ that coincides
with the initial function $c$ on $(0,t_\infty)$. Setting $t_M^c =
\inf\{ t>0 \vert c(t) = \gamma_M \}$ we see that $c(t_M^c)
= \gamma_M \downarrow-\infty$ as $M \rightarrow\infty$. Hence for
any $A \in{\mathcal B}(\mathbb{R})$ we see that $\mathsf P(B_{\tau
_{b,c}} \in A \cap
(c(t_M^c),\infty)) = \mathsf P(B_{\tau_{b,c_M}} \in A \cap
(c(t_M^c),\infty)) = \mu(A \cap(c(t_M^c),\infty)) \rightarrow
\mu(A)$ as $M \rightarrow\infty$ from where it follows that
$\mathsf P(B_{\tau_{b,c}} \in A) = \mu(A)$. This shows that
$B_{\tau_{b,c}} \sim\mu$ as required. Moreover we claim that this
is the only case we need to consider since if both $\{ t_N^b
\vert N \ge1 \}$ and $\{ t_N^c \vert N \ge1 \}$ are
bounded, then without loss of generality we can assume (by passing to
a subsequence if needed) that $t_N^c \rightarrow t_\infty^c <
\infty$ with $t_\infty^c > t_\infty$ first so that $t_N^b \uparrow
t_\infty$ as $N \rightarrow\infty$. In this case we can repeat the
preceding construction and extend $c$ to $[t_\infty,t_\infty^c)$ so
that we again have $B_{\tau_{b,c}} \sim\mu$ by the same argument.
If $t_\infty^c = t_\infty$, however, then the same argument as in the
case of $t_\infty= \infty$ above shows that the latter relation
also holds. Thus in both cases we have $t_N^b \le T$ and $t_N^c \le
T$ for all $N \ge1$ with $T:= t_\infty^c$ so that
$\mu((\gamma_N,\beta_N)) = \mathsf P(B_{\tau_{b,c}} \in
(\gamma_N,\beta_N)) = \mathsf P(B_{\tau_{b,c}} \in(c(t_N^c),b(t_N^b)))
\le\mathsf P(\tau_{b,c} \le T)$ for all $N \ge1$. Letting $N
\rightarrow
\infty$ and using that $\mu((\gamma_N,\beta_N)) \rightarrow1$, we
get $\mathsf P(\tau_{b,c} \le T) = 1$ which clearly is impossible since
$T<\infty$. It follows therefore that $B_{\tau_{b,c}} \sim\mu$ in
all possible cases and the proof is complete.\quad\qed
\end{longlist}\noqed
\end{pf}

%re2 #&#
\begin{rem}\label{rem2}
Note that $b$ from (I) and $c$ from (II) are
always finite valued since otherwise $\mu(\mathbb{R}_+)<1$ or $\mu
(\mathbb{R}_-)<1$,
respectively. Note also that either $b$ or $c$ from (III) can
formally take value $+\infty$ or $-\infty$, respectively, from some
time onwards; however, when this happens to either function, then the
other function must remain finite valued [note that (I) and (II) can
be seen as special cases of (III) in this sense too]. Note finally
that the result and proof of Theorem~\ref{teo1} including the same remarks
remain valid if $B_0 \sim\nu$ where $\nu$ is a probability measure
on $\mathbb{R}$ such that $\operatorname{supp}(\nu) \subseteq
[-p,q]$ with $\mu([-p,q])=0$
for some $p>0$ and $q>0$.
\end{rem}

%re3 #&#
\begin{rem}\label{rem3}
Since the arguments in the proof of Theorem~\ref{teo1} can
be repeated over any subsequence of $\{ b_n \vert n \ge1
\}$ or $\{ c_n \vert n \ge1 \}$ [when constructed with no
upper or lower bound on the partitions of $\operatorname{supp}(\mu)$
as well] it
follows that $B_{\tau_{b_n,c_n}}$ not only converges to
$B_{\tau_{b,c}}$ over a subsequence $\mathsf P$-a.s., but this
convergence also holds for the entire sequence in
$\mathsf P$-probability. Indeed, if this would not be the case, then for
some subsequence no further subsequence would converge $\mathsf P$-a.s.
The initial argument of this remark combined with the uniqueness
result of Theorem~\ref{teo10} below would then yield a contradiction. The
fact that $B_{\tau_{b_n,c_n}}$ always converges to $B_{\tau_{b,c}}$
in $\mathsf P$-probability as $n \rightarrow\infty$ makes the
derivation fully constructive and amenable to algorithmic
calculations described next.
\end{rem}

%re4 #&#
\begin{rem}\label{rem4}
The construction presented in the proof above
yields a simple algorithm for computing $b_n$ and $c_n$, which in
turn provide numerical approximations of $b$ and $c$. Key elements
of the algorithm can be described as follows. Below we let
$\varphi(x) = (1/\sqrt{2\pi}) e^{-x^2 /2}$ and $\Phi(x) =
(1/\sqrt{2\pi}) \int_{-\infty}^x e^{-y^2 /2} \,dy$ for $x \in
\mathbb{R}$ denote the standard normal density and distribution function,
respectively.

In the one-sided case (I) when $\operatorname{supp}(\mu) \subseteq
\mathbb{R}_+$ recall
the well-known expressions (cf. \cite{BS})
%
%
%e2.22 #&#
%e2.23 #&#
\begin{eqnarray}
\mathsf P(B_t \in dx, \tau_y>t) &=&
\frac{1}{\sqrt{t}} \biggl[ \varphi\biggl( \frac{x}{\sqrt{t}} \biggr) -
\varphi \biggl(\frac{x - 2y}{\sqrt{t}} \biggr) \biggr] \,dx
\nonumber\\[-8pt]\label{222}  \\[-8pt]
&=:& f(t,x,y) \,dx,\nonumber
\\
\label{223}  \mathsf P(\tau_y \le t) &=& 2 \biggl[ 1 - \Phi
\biggl( \frac{y}{\sqrt{t}} \biggr) \biggr] =: g(t,y)
\end{eqnarray}
for $t>0$ and $x<y$ with $y>0$ where we set $\tau_y = \inf\{
t>0 \vert B_t=y \}$. Using stationary and independent
increments of $B$ (its Markov property), we then read from Part~1.1
of the proof above that the algorithm runs as follows:
%
%
%e2.24 #&#
%e2.25 #&#
%e2.26 #&#
\begin{eqnarray}
\label{224} g_k(t) &:=& \int_{-\infty}^{x_{k-1}^n}
g\bigl(t,x_k^n - y\bigr) f_{k-1}(y) \,dy,
\\
\label{225} t_k^n &:=& t_{k-1}^n
+ \inf\bigl\{ t>0 \vert g_k(t) = p_k^n
\bigr\},
\\
\label{226} f_k(x) &:=& \int_{-\infty}^
{x_{k-1}^n}
f\bigl(t_k^n - t_{k-1}^n,x - y,
x_k^n - y\bigr) f_{k-1}(y) \,dy
\end{eqnarray}
for $k = 1,2, \ldots, m_n$ where we initially set $t_0:=0$, $x_0:=0$
and $f_0(x) \,dx:= \delta_0(dx)$. This yields the time points
$t_1^n, t_2^n, \ldots, t_{m_n}^n$ which determine $b_n$ by the
formula
%
%
%e2.27 #&#
\begin{equation}
\label{227} b_n(t) = \sum_{k=1}^{m_n}
x_k^n I\bigl(t_{k-1}^n < t \le
t_k^n\bigr)
\end{equation}
for $t \ge0$. The algorithm is stable and completes within a
reasonable time frame; see Figure~\ref{fig2} below for the numerical output
when the target law $\mu$ is exponentially distributed with
intensity 1.

%%%%%%%%%%%%%%%%%%%%%%%%%%%%%%%%%%%%%%%%%%%%%%%%%%%%%%%%%%%%%%%%%%%%%%%%%%%%%%%%%
%%% Figure 2 %%%
%%%%%%%%%%%%%%%%%%%%%%%%%%%%%%%%%%%%%%%%%%%%%%%%%%%%%%%%%%%%%%%%%%%%%%%%%%%%%%%%%

%
%f2 #&#
\begin{figure}%[t]

\includegraphics{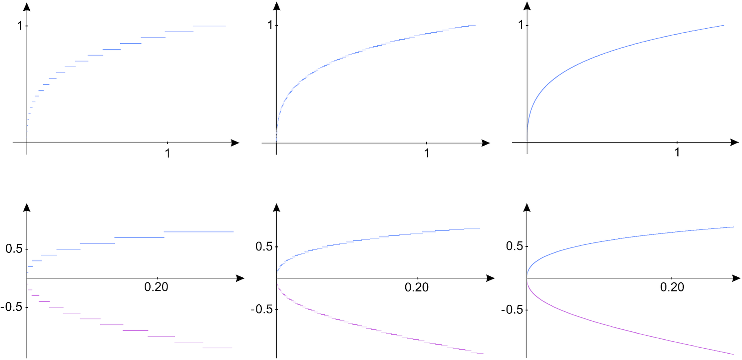}

\caption{Functions $b_n$ and $c_n$ calculated using the
algorithm from the proof of Theorem~\protect\ref{teo1} as described in Remark~\protect\ref{rem4}. The
first row corresponds to the target law $\mu$ which is exponentially
distributed with intensity $1$ for $n ={}$20, 100, 500,
respectively, with equidistant partition of $\mathbb{R}_+$ having the step
size equal to $1/n$ and the number of time points $m_n$ equal to
$n$. The second row corresponds to the target law $\mu$ which is
normally distributed with mean $1$ and variance $1$ for $n ={}$10,
50, 250, respectively, with equidistant partition of $\mathbb{R}$ having
the step size equal to $1/n$ and the number of time points $m_n +
l_n$ equal to $2n$.}\vspace*{-3pt}\label{fig2}
\end{figure}

%%%%%%%%%%%%%%%%%%%%%%%%%%%%%%%%%%%%%%%%%%%%%%%%%%%%%%%%%%%%%%%%%%%%%%%%%%%%%%%%%
%%% End Figure %%%
%%%%%%%%%%%%%%%%%%%%%%%%%%%%%%%%%%%%%%%%%%%%%%%%%%%%%%%%%%%%%%%%%%%%%%%%%%%%%%%%%

In the two-sided case (III) when $\operatorname{supp}(\mu) \subseteq
\mathbb{R}$ recall
the well-known expressions (cf. \cite{BS})
%
%
%e2.28 #&#
%e2.29 #&#
%e2.30 #&#
\begin{eqnarray}
\qquad && \mathsf P(B_t \in dx, \tau_{y,z}>t)\nonumber
\\
\label{228}  &&\qquad = \frac{1}{\sqrt{t}} \sum_{n=-
\infty}^\infty\biggl[
\varphi\biggl( \frac{x + 2n(y - z)}{\sqrt{t}} \biggr) - \varphi\biggl
(\frac{x + 2n(y - z) - 2y}{\sqrt{t}}
\biggr) \biggr] \,dx
\\
&&\qquad  =: f(t,x,y,z) \,dx,\nonumber
\\
&& \mathsf P(\tau_y < \tau_z, \tau_{y,z} \le t)\nonumber
\\
\label{229} &&\qquad = 2 \sum_{n=0}^\infty
\biggl[ \Phi\biggl( \frac{(2n + 1)(y - z) -
z}{\sqrt{t}} \biggr) - \Phi\biggl(\frac{(2n + 1)(y - z) + z}{
\sqrt{t}}\biggr) \biggr]
\\
&&\qquad  =: g(t,y,z),\nonumber
\\
&& \mathsf P(\tau_z < \tau_y, \tau_{y,z} \le t)\nonumber
\\
\label{230} &&\qquad = 2 \sum_{n=0}^\infty
\biggl[ \Phi\biggl( \frac{(2n + 1)(y - z)
+ y}{\sqrt{t}} \biggr) - \Phi\biggl(\frac{(2n + 1)(y - z) - y}{
\sqrt{t}} \biggr) \biggr]
\\
&&\qquad  =: h(t,y,z)\nonumber
\end{eqnarray}
for $t>0$ and $z < x < y$ with $z < 0 < y$ where we set $\tau_w =
\inf\{ t>0 \vert B_t=w \}$ for $w \in\{y,z \}$ and
$\tau_{y,z} = \tau_y \wedge\tau_z$. Using stationary and
independent increments of~$B$ (its Markov property), we then read
from Part~3.1 of the proof above that the algorithm runs as follows:
%
%
%e2.31 #&#
%e2.32 #&#
%e2.33 #&#
%e2.34 #&#
\begin{eqnarray}
\label{231} \qquad g_k(t)&:=& \int_{\bar y_{k-1}^n}^{\bar x_{k-1}^n}
g\bigl(t,\bar x_k^n - z,\bar y_k^n
- z\bigr) f_{k-1}(z) \,dz,
\\
\label{232} h_k(t)&:=& \int_{\bar y_{k-1}^n}^{\bar x_{k-1}^n}
h\bigl(t,\bar x_k^n - z,\bar y_k^n
- z\bigr) f_{k-1}(z) \,dz,
\\
\label{233}  t_k^n&:=& t_{k-1}^n
+ \bigl( \inf\bigl\{ t>0 \vert g_k(t) = \bar p_k^n
\bigr\} \wedge\inf\bigl\{ t>0 \vert h_k(t) = \bar
q_k^n \bigr\} \bigr),
\\
\label{234} f_k(x)&:=& \int_{\bar y_{k-1}^n}^{\bar x_{k-1}^n}
f\bigl(t_k^n - t_{k-1}^n,x - z,\bar
x_k^n - z,\bar y_k^n - z\bigr)
f_{k-1}(z) \,dz
\end{eqnarray}
for $k = 1,2, \ldots, m_n + l_n$ where we initially set $t_0:= 0$,
$\bar x_0^n:= 0$, $\bar y_0^n:= 0$, $\bar x_1^n:= x_1^n$, $\bar
y_1^n:= y_1^n$, $f_0(x) \,dx:= \delta_0(dx)$ and denoting the
first infimum in (\ref{233}) by $I_k^n$ and the second infimum in
(\ref{233}) by $J_k^n$, this is\vspace*{1pt} then continued as follows: if
$I_k^n>J_k^n$, then $\bar x_{k+1}^n:= \inf\{ x_l \vert x_l
> \bar x_k^n \}$, $\bar y_{k+1}^n:= \bar y_k^n$,\vspace*{1pt} $\bar p_{k+1}^n:=
p(\bar x_{k+1}^n)$, $\bar q_{k+1}^n:= \bar q_k^n -
h_k(I_k^n)$; if $J_k^n>I_k^n$,\vspace*{1pt} then $\bar y_{k+1}^n:= \sup\{
y_l \vert y_l < \bar y_k^n \}$, $\bar x_{k+1}^n:= \bar x_k^n$, $\bar
q_{k+1}^n:= q(\bar y_{k+1}^n)$,
$\bar p_{k+1}^n:= \bar p_k^n - g_k(J_k^n)$; if $I_k^n=J_k^n$,\vspace*{1pt} then
$\bar x_{k+1}^n:= \inf\{ x_l \vert x_l > \bar x_k^n \}$, $\bar
y_{k+1}^n:= \sup\{ y_l \vert y_l < \bar y_k^n \}$,
$\bar p_{k+1}^n:= p(\bar x_{k+1}^n)$, $\bar q_{k+1}^n:= q(\bar
y_{k+1}^n)$ where\vspace*{1pt} we set $p(x) = p_k^n$ for $x=x_k^n$ and \mbox{$q(y) =
q_k^n$} for $y=y_k^n$. This yields the time points $t_1^n, t_2^n,
\ldots, t_{m_n+l_n}^n$ which determine $b_n$ and $c_n$ by the
formulae
%
%
%e2.35 #&#
\begin{eqnarray}\label{235}
b_n(t) &=& \sum_{k=1}^{m_n+l_n}
\bar x_k^n I\bigl(t_{k-1}^n < t \le
t_k^n\bigr) \quad\mbox{and}
\nonumber\\[-8pt]\\[-8pt]
c_n(t) &=&
\sum_{k=1}^{m_n+l_n} \bar y_k^n
I\bigl(t_{k-1}^n < t \le t_k^n
\bigr)\nonumber
\end{eqnarray}
for $t \ge0$. The algorithm is stable and completes within a
reasonable time frame; see Figure~\ref{fig2} above for the numerical
output when the target law $\mu$ is normally distributed with mean
$1$ and variance $1$.
\end{rem}

%re5 #&#
\begin{rem}\label{rem5}
Note that $\tau_b$ from (I) could also be defined
by
%
%
%e2.36 #&#
\begin{equation}
\label{236} \tau_b = \inf\bigl\{ t>0 \vert B_t =
b(t) \bigr\}
\end{equation}
and that $B_{\tau_b} = b(\tau_b)$. This is easily verified since $b$
is left-continuous and increasing. The same remark applies to
$\tau_c$ from (II) and $\tau_{b,c}$ from (III) with $B_{\tau_{b,c}}$
being equal to $b(\tau_{b,c})$ or $c(\tau_{b,c})$. From (\ref{28})
and (\ref{29}) we also see that these inequalities and equalities
in the definitions of the stopping times can be replaced by strict
inequalities and that all relations remain valid almost surely in
this case. Similarly, in all these definitions we could replace
left-continuous functions $b$ and $c$ with their right-continuous
versions defined by $b(t):= b(t+)$ and $c(t):= c(t+)$ for $t>0$,
respectively. All previous facts in this remark remain valid in this
case too.
\end{rem}

%re6 #&#
\begin{rem}\label{rem6}
If $\mu(\{ 0 \}) =: p > 0$ in Theorem~\ref{teo1}, then we
can generate a random variable $\zeta$ independently from $B$ such
that $\zeta$ takes two values $0$ and $\infty$ with probabilities
$p$ and $1 - p$, respectively. Performing the same construction with
the stopped sample path $t \mapsto B_{t \wedge\zeta}$ yields the
existence of functions $b$ and $c$ as in Theorem~\ref{teo1} with $B^\zeta=
(B_{t \wedge\zeta})_{t \ge0}$ in place of $B=(B_{t \ge0})_{t \ge
0}$. The resulting stopping time may be viewed as randomised through
the initial condition.
\end{rem}

%re7 #&#
\begin{rem}\label{rem7}
Two main ingredients in the proof of Theorem~\ref{teo1}
above are (i)~embedding in discrete laws and (ii) passage to the
limit from discrete to general laws. If the standard Brownian motion
$B$ is replaced by a continuous (time-homogeneous) Markov process
$X$, we see from the proof above that (i) can be achieved when
%
%
%e2.37 #&#
\begin{equation}
\label{237} t \mapsto\mathsf P_{ x}(\tau_y <
\tau_z, \tau_{y,z} \le t) \quad\mbox{and}\quad t \mapsto
\mathsf P_{ x}(\tau_z < \tau_y,
\tau_{y,z} \le t)
\end{equation}
are continuous on $\mathbb{R}_+$ and $\mathsf P_{ x}(\tau_{y,z} > t)
\downarrow
0$ as $t \uparrow\infty$ for all $-\infty\le z < x < y \le\infty$
with $\vert z \vert\wedge\vert y \vert< \infty$ and
$\mathsf P_{ x}(X_0=x)=1$ where we set $\tau_w = \inf\{ t>0
\vert X_t=w \}$ for $w \in\{y,z \}$ and $\tau_{y,z} = \tau_y
\wedge\tau_z$. We also see from the proof above that (ii) can be
achieved when
%
%
%e2.38 #&#
\begin{equation}
\label{238} \tau_b = \sigma_b \quad\mathsf P_{ 0}\mbox{-a.s.} \quad\mbox{and}\quad\tau_c =
\sigma_c\quad\mathsf P_{ 0}\mbox{-a.s.},
\end{equation}
where the first equality holds for any left-continuous increasing
function $b$ with $\tau_b = \inf\{ t>0 \vert X_t \ge b(t)
\}$ and $\sigma_b = \inf\{ t>0 \vert X_t > b(t) \}$, and
the second equality holds for any left-continuous decreasing
function $c$ with $\tau_c = \inf\{ t>0 \vert X_t \le c(t)
\}$ and $\sigma_c = \inf\{ t>0 \vert X_t < c(t) \}$. In
particular, by verifying (\ref{237}) and (\ref{238}) in the proof
of Corollary~\ref{cor8} below we will establish that the result of Theorem~\ref{teo1}
extends to all recurrent diffusion processes $X$ in the sense of
It\^o and McKean \cite{IM} (see \cite{BS}, Chapter II, for a review).
This extension should also hold for nonrecurrent diffusion
processes $X$ and ``admissible'' target laws $\mu$ (cf. \cite{PP}) as
well as for more general standard Markov processes $X$ satisfying
suitable modifications of (\ref{237}) and (\ref{238}) in the
admissible setting. We leave precise formulations of these more
general statements and proofs as informal conjectures open for
future developments.
\end{rem}

%co8 #&#
\begin{cor}\label{cor8}
The result of Theorem~\ref{teo1} remains valid if
the standard Brownian motion $B$ is replaced by any recurrent
diffusion process $X$.
\end{cor}

\begin{pf} As pointed out above the proof can be carried out in
the same way as the proof of Theorem~\ref{teo1} if we show that (\ref{237})
and (\ref{238}) are satisfied. Note that $\mathsf P_{ x}(\tau_{y,z}
> t) \downarrow0$ as $t \uparrow\infty$ for all $-\infty\le z < x
< y \le\infty$ with $\vert z \vert\wedge\vert y \vert< \infty$
since $X$ is recurrent. Recall also that all recurrent diffusions
are regular (see \cite{BS}, Chapter~II, for definitions).
\begin{longlist}[1.]
\item[1.] We first show that the functions in (\ref{237}) are continuous.
Clearly by symmetry it is enough to show that the first function is
continuous. For this, set $F(t) = \mathsf P_{ x}(\tau_y < \tau_z,
\tau_{y,z} \le t)$ for $t \ge0$ where $-\infty\le z < x < y <
\infty$ are given and fixed. Since $t \mapsto F(t)$ is increasing
and right-continuous we see that it is enough to disprove the
existence of $t_1>0$ such that $F(t_1) - F(t_1-) = \mathsf P_{ x}(\tau_y
< \tau_z, \tau_y = t_1) > 0$. Since this implies that
$\mathsf P_{ x}(\tau_y = t_1) > 0$ we see that it is enough to show that
the distribution function $t \mapsto\mathsf P_{ x}(\tau_y \le t)$ is
continuous for $x<y$ in $\mathbb{R}$ given and fixed. For this, let $p$
denote the transition density of $X$ with respect to its speed
measure $m$ in the sense that $\mathsf P_{ x}(X_t \in A) = \int_A
p(t;x,y) m(dy)$ holds for all $t>0$ and all $A \in{\mathcal B}(\mathbb
{R})$. It is
well known (cf. \cite{IM}, page 149) that $p$ may be chosen to be
jointly continuous (in all three variables). Next note that for any
$s>0$ given and fixed the mapping $t \mapsto\mathsf E _x [
\mathsf P_{ X_s}(\tau_y \le t) ] = \int_{\mathbb{R}} \mathsf P_{
z}(\tau_y \le t)
p(s;x,z) m(dz)$ is increasing and right-continuous on $(0,\infty)$
so that $G(t,s):= \mathsf E _x [ \mathsf P_{ X_s}(\tau_y = t)
] =
\int_{\mathbb{R}} \mathsf P_{ z}(\tau_y = t) p(s;x,z) m(dz) = 0$ for
all $t
\in(0,\infty) \setminus C_s$ where the set $C_s$ is at most
countable. Setting $C:= \bigcup_{s \in\mathbb{Q}_+} C_s$ where
$\mathbb{Q}_+$ denotes the set of rational numbers in $(0,\infty)$,
we see that the set $C$ is at most countable and $G(t,s)=0$ for all
$t \in(0,\infty) \setminus C$ and all $s \in\mathbb{Q}_+$. Since
each $z \mapsto p(s;x,z)$ is a density function integrating to $1$
over $m(dz)$, and $s \mapsto p(s;x,z)$ is continuous on $(0,\infty)$,
we see by Scheff\'e's theorem (see, e.g., \cite{Bi}, page 215) that
$G(t,s_n) \rightarrow G(t,s)$ as $s_n \rightarrow s$ in $(0,\infty)$
for any $t>0$ fixed. Choosing these $s_n$ from $\mathbb{Q}_+$ for
given $s>0$ it follows therefore that $G(t,s)=0$ for all $t \in
(0,\infty) \setminus C$ and all $s>0$. By the Markov property we
moreover see that $\mathsf P_{ x}(\tau_y=t + s) \le\mathsf P_{
x}(\tau_y
\circ\theta_s = t) = G(t,s) = 0$ and hence $\mathsf P_{ x}(\tau_y=t +
s)=0$ for all $t \in(0,\infty) \setminus C$ and all $s>0$. Since
the set $C$ is at most countable it follows that
$\mathsf P_{ x}(\tau_y=t)=0$ for all $t>0$. This implies that $F$ is
continuous, and the proof of (\ref{237}) is complete.
\end{longlist}

\begin{longlist}[2.]
\item[2.] We next show that the equalities in (\ref{238}) are satisfied.
Clearly by symmetry it is enough to derive the first equality. Note
that replacing $B$ by $X$ in the proof of~(\ref{29}) above and
using exactly the same arguments yields the first equality in~(\ref{238}), provided that (\ref{214}) is established for $X$ in
place of $B$. This shows that the first equality in (\ref{238})
reduces to establishing that
%
%
%e2.39 #&#
\begin{equation}
\label{239} \mathsf P_{ 0}\bigl(\sigma_b^\delta
> t\bigr) \le\mathsf P_{ 0}\bigl( \tau_b^\delta>
t\bigr)
\end{equation}
for all $t>0$ where $\sigma_b^\delta= \inf\{ t > \delta
\vert X_t > b(t) \}$ and $\tau_b^\delta= \inf\{ t >
\delta\vert X_t \ge b(t) \}$ for $\delta>0$ given and fixed.
Observe that $\sigma_b^\delta$ coincides with $\tau_{b+}^\delta:=
\lim_{ \varepsilon\downarrow0} \tau_{b+\varepsilon}^\delta$ where
$\tau_{b+\varepsilon}^\delta= \inf\{ t > \delta\vert  X_t
\ge
b(t) + \varepsilon\}$ as is easily seen from the definitions so that
(\ref{239}) is indeed equivalent to (\ref{214}) as stated above.

To establish (\ref{239}) consider first the case when $b$ is flat
on some time interval $I \subseteq(\delta,\infty)$, and denote the
joint value of $b$ on $I$ by $y$ meaning that $b(t)=y$ for all $t
\in I$. Consider the stopping times $\tau:= \inf\{ t>\delta
\vert X_t = y \}$ and $\sigma:= \inf\{ t>\tau\vert
X_t > y \}$. Since $X$ is recurrent we know that both $\tau$ and
$\sigma$ are finite valued under $\mathsf P_{ 0}$. Note that $\sigma=
\tau+ \rho\circ\theta_\tau$ where $\rho:= \inf\{ t>0
\vert X_t > y \}$ is a stopping time. By the strong Markov
property of $X$ applied at $\tau$, we thus have $\mathsf P_{ 0}(
\sigma
= \tau) = \mathsf P_{ 0}( \rho\circ\theta_\tau= 0) =
\mathsf P_{ X_\tau}( \rho= 0) = \mathsf P_{ y}( \rho= 0) = 1$ where
the final equality follows since $X$ is regular (cf. \cite{BS}, page~13). Hence we see that $X_{\tau+t} > y$ for infinitely many $t$
in each $(0,\varepsilon]$ for $\varepsilon>0$ with $\mathsf P_{ 0}$-probability one.
In particular, this shows that on the set $\{ \sigma_b^\delta
> t \}$ with $t>0$ given and fixed the sample path of $X$ stays
strictly below $b$ on the time interval $I \setminus\sup(I)$ with
$\mathsf P_{ 0} $-probability one for each time interval $I \subseteq
(\delta,t)$ on which $b$ is flat. Since $(\delta,t)$ can be written
as a countable union of disjoint intervals on each of which $b$ is
either flat or strictly increasing, we see that the previous
conclusion implies that
%
%
%e2.40 #&#
\begin{eqnarray}
\label{240} \mathsf P_{ 0}\bigl(\sigma_b^\delta
> t\bigr) &\le&\mathsf P_{ 0}\bigl( X_s < b(s + h)
\mbox{ for all } s \in(\delta,t) \bigr)\nonumber
\\
&\le&\mathsf P_{ 0}\bigl( X_{r-h} < b(r)
\mbox{ for all } r \in(\delta+ h,t + h) \bigr)
\\
\nonumber
&\le&\mathsf P_{ 0}\bigl( X_{r-h} < b(r) \mbox{ for
all } r \in(\delta+ h_0,t] \bigr)
\end{eqnarray}
for any $h \in(0,h_0)$ where $h_0 \in(0,\delta/2)$ is given and
fixed. By the Markov property and Scheff\'e's theorem applied as
above, we find that
%
%
%e2.41 #&#
\begin{eqnarray}
\label{241} \qquad && \mathsf P_{ 0}\bigl( X_{r-h} < b(r)
\mbox{ for all } r \in(\delta+ h_0,t] \bigr)\nonumber
\\
&&\qquad  = \mathsf E _0 \bigl[ \mathsf P_{
X_{\delta/2
-h}}\bigl(
X_{r-\delta/2} < b(r) \mbox{ for all } r \in (\delta+
h_0,t] \bigr) \bigr]\nonumber
\\
&&\qquad  = \int_{\mathbb{R}} \mathsf P_{ y}\bigl(
X_{r-\delta/2} < b(r) \mbox{ for all } r \in(\delta+ h_0,
t] \bigr) p(\delta/2 - h;0,y) m(dy)
\nonumber\\[-8pt]\\[-8pt]
&&\qquad  \longrightarrow\int_{\mathbb{R}} \mathsf
P_{ y}\bigl( X_{r-\delta/2} < b(r) \mbox{ for all } r \in(\delta+
h_0,t] \bigr) p(\delta/2;0,y) m(dy)\nonumber
\\
&&\qquad  = \mathsf E _0 \bigl[ \mathsf P_{ X_{\delta/2}} \bigl(
X_{r-\delta/2} < b(r) \mbox{ for all } r \in(\delta+
h_0,t] \bigr) \bigr]\nonumber
\\
&&\qquad  =\mathsf P_{ 0}\bigl( X_r < b(r)\mbox{ for
all } r \in(\delta+ h_0,t] \bigr)\nonumber
\end{eqnarray}
as $h \downarrow0$. Combining (\ref{240}) and (\ref{241}) we get
%
%
%e2.42 #&#
\begin{equation}
\label{242} \mathsf P_{ 0}\bigl(\sigma_b^\delta
> t\bigr) \le\mathsf P_{ 0}\bigl( X_r < b(r) \mbox{ for
all } r \in(\delta+ h_0,t] \bigr)
\end{equation}
for all $h_0 \in(0,\delta/2)$. Letting $h_0 \downarrow0$ in
(\ref{242}) we find that
%
%
%e2.43 #&#
\begin{equation}
\label{243} \mathsf P_{ 0}\bigl(\sigma_b^\delta
> t\bigr) \le\mathsf P_{ 0}\bigl( X_r < b(r) \mbox{ for
all } r \in(\delta,t] \bigr) = \mathsf P_{ 0}\bigl( \tau
_b^\delta> t\bigr)
\end{equation}
for all $t>0$. This establishes (\ref{239}) and hence $\tau_b =
\sigma_b$ $\mathsf P_{ 0} $-a.s. as explained above. The proof of
(\ref{238}) is therefore complete.\quad\qed
\end{longlist}\noqed
\end{pf}

Note that the claims of Remarks~\ref{rem2}--\ref{rem6} extend to the setting of
Corollary~\ref{cor8} with suitable modifications in Remark~\ref{rem4} since the
process no longer has stationary and independent increments and some
of the expressions may no longer be available in closed form.

In the setting of Theorem~\ref{teo1} or Corollary~\ref{cor8}, let $F_\mu$ denote the
distribution function of $\mu$. The following proposition shows that
(i) jumps of $b$ or $c$ correspond exactly to flat intervals of
$F_\mu$ (i.e., no mass of $\mu$), and (ii) flat intervals of $b$ or
$c$ correspond exactly to jumps of $F_\mu$ (i.e., atoms of $\mu$).
In particular, from (i) we see that if $F_\mu$ is strictly
increasing on $\mathbb{R}_+$, then $b$ is continuous, and if $F_\mu$ is
strictly increasing on $\mathbb{R}_-$, then $c$ is continuous. Similarly,
from (ii) we see that if $F_\mu$ is continuous on $\mathbb{R}_+$,
then $b$ is
strictly increasing, and if $F_\mu$ is continuous on $\mathbb{R}_-$,
then $c$
is strictly decreasing.

%pr9 #&#
\begin{prop}[(Continuity)]\label{prop9} In the setting of Theorem~\ref{teo1} or Corollary~\ref{cor8} we have:
%
%
%e2.44 #&#
%e2.45 #&#
%e2.46 #&#
%e2.47 #&#
\begin{eqnarray}
\label{244} \qquad && b(t+) > b(t)\quad\mbox{if and only if}\quad\mu\bigl(
\bigl(b(t),b(t+)\bigr)\bigr) = 0,
\\
\label{245} && b(t) = b(t - \varepsilon)\quad\mbox{for some }
\varepsilon>0
\quad\mbox{if and only if}\quad\mu\bigl(\bigl\{ b(t) \bigr\}\bigr) > 0,
\\
\label{246} && c(t+) < c(t) \quad\mbox{if and only if}\quad\mu\bigl(
\bigl(c(t+),c(t)\bigr)\bigr) = 0,
\\
\label{247} && c(t) = c(t - \varepsilon) \quad\mbox{for some }
\varepsilon
>0 \quad\mbox{if and only if}\quad\mu\bigl(\bigl\{ c(t) \bigr\}\bigr)
> 0,
\end{eqnarray}
for any $t>0$ given and fixed.
\end{prop}

\begin{pf} All statements follow from the construction and
basic properties of $b$ and $c$ derived in the proof of Theorem~\ref{teo1}.
\end{pf}

%s3 #&#
\section{Uniqueness}\label{sec3}
%%%%%%%%%%%%%%%%%%%%

In this section we state and prove the main uniqueness result. Note
that the result and proof remain valid in the more general case
addressed at the end of Remark~\ref{rem2}, and the method of proof is also
applicable to more general processes (cf. Remark~\ref{rem7}).

%th10 #&#
\begin{teo}[(Uniqueness)]\label{teo10}
 In the setting of Theorem~\ref{teo1}
or Corollary~\ref{cor8} the functions $b$ and $c$ are uniquely determined by
the law $\mu$.
\end{teo}

\begin{pf}
To simplify the exposition we will derive (I) in
full detail. It is clear from the proof below that the same
arguments can be used to derive (II) and (III).
\begin{longlist}[1.]
\item[1.] Let us assume that $b_1\dvtx (0,\infty) \rightarrow\mathbb{R}_+$
and $b_2\dvtx
(0,\infty) \rightarrow\mathbb{R}_+$ are left-continuous increasing
functions such that $X_{\tau_{b_1}} \sim\mu$ and $X_{\tau_{b_2}}
\sim\mu$ where $\tau_{b_1} = \inf\{ t>0 \vert X_t \ge
b_1(t) \}$ and $\tau_{b_2} = \inf\{ t>0 \vert X_t \ge
b_2(t) \}$. We then need to show that $b_1 = b_2$. For this, we
will first show that $b:= b_1 \wedge b_2$ also solves the embedding
problem in the sense that $X_{\tau_b} \sim\mu$ where $\tau_b =
\inf\{ t>0 \vert X_t \ge b(t) \}$. The proof of this fact
can be carried out as follows.
\end{longlist}

\begin{longlist}[2.]
\item[2.] Let $A = \{ x \in\operatorname{supp}(\mu) \vert\mu(\{x\})>0 \}
$ and for any given $x \in A$ set $\ell_i(x) = \inf\{ t \in
(0,\infty) \vert b_1(t)=x \}$ and $r_i(x) = \sup\{ t \in
(0,\infty) \vert b_1(t)=x \}$ when $i=1,2$. By (\ref{245}) we
know that $[\ell_i(x),r_i(x)]$ is a nonempty interval. Moreover,
note that the functions $\ell_i$ and $r_i$ are also well defined on
$\operatorname{supp}(\mu) \setminus A$ (with the convention $\inf
\varnothing= \sup\varnothing= +\infty$) in which case we have
$\ell_i = r_i$ for $i=1,2$. With this notation in mind consider the
sets
%
%
%e3.1 #&#
%e3.2 #&#
%e3.3 #&#
%e3.4 #&#
%e3.5 #&#
%e3.6 #&#
%e3.7 #&#
\begin{eqnarray}
\label{31} G_{1,1} &=& \bigl\{ x \in\operatorname{supp}(\mu)
\setminus A \vert\ell_1(x) < \ell_2(x) \bigr\},
\\
G_{1,2} &=& \bigl\{ x \in A \vert r_1(x) <r_2(x) \bigr\},
\\
G_{1,3} &=& \bigl\{ x \in A \vert\ell_1(x) <\ell_2(x) \mbox{ and } r_1(x)=r_2(x)\bigr\},
\\
G_{2,1} &=& \bigl\{ x \in\operatorname{supp}(\mu) \setminus A \vert \ell_2(x) \le\ell_1(x) \bigr\},
\\
G_{2,2} &=& \bigl\{ x \in A \vert r_2(x) <r_1(x) \bigr\},
\\
G_{2,3} &=& \bigl\{ x \in A \vert\ell_2(x) <\ell_1(x) \mbox{ and } r_1(x)=r_2(x)\bigr\},
\\
G_{2,4} &=& \bigl\{ x \in A \vert\ell_1(x) =\ell_2(x) \mbox{ and } r_1(x)=r_2(x)\bigr\}.
\end{eqnarray}
Set $G_1:= G_{1,1} \cup G_{1,2} \cup G_{1,3}$ and $G_2:= G_{2,1}
\cup G_{2,2} \cup G_{2,3} \cup G_{2,4}$. Note that $G_1$ and $G_2$
are disjoint and $\operatorname{supp}(\mu) = G_1 \cup G_2$. Setting
$\tau_1:=
\tau_{b_1}$ and $\tau_2:= \tau_{b_2}$ we claim that
%
%
%e3.8 #&#
\begin{equation}
\label{38} \mathsf P(X_{\tau_1} \in G_1,
X_{\tau_2} \in G_2) = 0.
\end{equation}
Indeed, if $X_{\tau_1} \in G_1$, then $X_{\tau_1} = b_1(\tau_1) \ge
b_2(\tau_1)$ so that $\tau_2 \le\tau_1$, while if $X_{\tau_2} \in
G_2$, then $X_{\tau_2} = b_2(\tau_2) \ge b_1(\tau_2)$ so that $\tau_1
\le\tau_2$. Since $G_1$ and $G_2$ are disjoint, this shows that the
set in (\ref{38}) is empty and thus has $\mathsf P$-probability zero as
claimed. From (\ref{38}) we see that
%
%
%e3.9 #&#
\begin{equation}
\label{39} \mathsf P(X_{\tau_1} \in G_1) = \mathsf
P(X_{\tau_1} \in G_1, X_{\tau_2} \in
G_1).
\end{equation}
Since $X_{\tau_1} \sim X_{\tau_2}$ this is further equal to
%
%
%e3.10 #&#
\begin{equation}
\label{310} \mathsf P(X_{\tau_2} \in G_1) = \mathsf
P(X_{\tau_2} \in G_1, X_{\tau_1} \in
G_1) + \mathsf P(X_{\tau_2} \in G_1,
X_{\tau_1} \in G_2)
\end{equation}
from where we also see that
%
%
%e3.11 #&#
\begin{equation}
\label{311} \mathsf P(X_{\tau_1} \in G_2,
X_{\tau_2} \in G_1) = 0.
\end{equation}
It follows therefore that
%
%
%e3.12 #&#
\begin{equation}
\label{312} \mathsf P(X_{\tau_1} \in G_2) = \mathsf
P(X_{\tau_1} \in G_2, X_{\tau_2} \in G_2).
\end{equation}
From (\ref{39}) and (\ref{312}) we see that the sets $\Omega_1 =
\{ X_{\tau_1} \in G_1, X_{\tau_2} \in G_1 \}$ and
$\Omega_2 = \{ X_{\tau_1} \in G_2, X_{\tau_2} \in G_2
\}$ form a partition of $\Omega$ with $\mathsf P$-probability one.
Moreover, note that for $\omega\in\Omega_1$ we have
$X_{\tau_1}(\omega) \in G_1$ so that $\tau_2(\omega) \le
\tau_1(\omega)$ and hence $\tau_b(\omega) = \tau_2(\omega)$, and for
$\omega\in\Omega_2$ we have $X_{\tau_2}(\omega) \in G_2$ so that
$\tau_1(\omega) \le\tau_2(\omega)$ and hence $\tau_b(\omega) =
\tau_1(\omega)$. This implies that for every $C \in{\mathcal
B}(\operatorname{supp}(\mu))$
we have
%
%
%e3.13 #&#
\begin{eqnarray}\label{313}
&& \mathsf P(X_{\tau_b} \in C)\nonumber
\\
&&\qquad  = \mathsf P\bigl( \{
X_{\tau_2} \in C \} \cap\Omega_1\bigr) + \mathsf P\bigl( \{
X_{\tau_1} \in C \} \cap\Omega_2\bigr)\nonumber
\\
&&\qquad = \mathsf P(X_{\tau_1} \in G_1,
X_{\tau_2} \in C \cap G_1) + \mathsf P(X_{\tau_1} \in
C \cap G_2, X_{\tau_2} \in G_2)
\\
&&\qquad  = \mathsf P(X_{\tau_2} \in C \cap G_1) + \mathsf
P(X_{\tau_1} \in C \cap G_2)\nonumber
\\
&&\qquad  = \mu(C \cap G_1) + \mu(C \cap G_2) = \mu(C),\nonumber
\end{eqnarray}
where we also use (\ref{311}) in the third equality. This shows
that $X_{\tau_b} \sim\mu$ as claimed.
\end{longlist}

\begin{longlist}[3.]
\item[3.] To complete the proof we can now proceed as follows. Since $b
\le b_i$ we see that $X_{\tau_b} \le X_{\tau_{b_i}}$ for $i=1,2$.
Moreover, since $X_{\tau_b} \sim X_{\tau_{b_i}}$ from the latter
inequality, we see that $X_{\tau_b} = X_{\tau_{b_i}}$ $\mathsf P$-a.s.
for $i=1,2$. As clearly this is not possible if for some $t>0$ we
would have $b_1(t) \ne b_2(t)$, it follows that $b_1=b_2$ and the
proof is complete.\quad\qed
\end{longlist}\noqed
\end{pf}

%s4 #&#
\section{Minimality}\label{sec4}
%%%%%%%%%%%%%%%%%%%%

In this section we show that the stopping time from Theorem~\ref{teo1} or
Corollary~\ref{cor8} is minimal in the sense of Monroe; see \cite{Mo}, page 1294.

%pr11 #&#
\begin{prop}[(Minimality)]\label{prop11}
In the setting of
Theorem~\ref{teo1} or Corollary~\ref{cor8} let $\tau= \tau_{b,c}$ with $c = -\infty$
if $\operatorname{supp}(\mu) \subseteq\mathbb{R}_+$ and $b =
+\infty$ if $\operatorname{supp}(\mu)
\subseteq\mathbb{R}_-$. Let $\sigma$ be any stopping time such that
%
%
%e4.1 #&#
%e4.2 #&#
\begin{eqnarray}
\label{41} X_\sigma&\sim& X_\tau,
\\
\label{42} \sigma&\le&\tau\quad\mathsf P\mbox{-a.s.}
\end{eqnarray}
Then $\sigma= \tau$ $\mathsf P$-a.s.
\end{prop}

\begin{pf}
Since $\int_0^N \mathsf P(\sigma\ge t) \,dt = \mathsf
E (\sigma
\wedge N) \le\mathsf E (\tau\wedge N) = \int_0^N \mathsf P(\tau\ge
t) \,dt$
for all $N \ge1$ by (\ref{42}) above, we see that it is enough to
show that $\mathsf P(\sigma\ge t) \ge\mathsf P(\tau\ge t)$ or equivalently
%
%
%e4.3 #&#
\begin{equation}
\label{43} \mathsf P(\sigma< t) \le\mathsf P(\tau< t)
\end{equation}
for all $t>0$. For this, note that from (\ref{41}) and (\ref{42})
combined with the facts that $b$ and $c$ are left-continuous
increasing and decreasing functions, respectively, it follows that
%
%
%e4.4 #&#
\begin{eqnarray}
\label{44} \mathsf P(\sigma< t) & =& \mathsf P \bigl( \sigma< t,
X_\sigma\in\bigl(c(t),b(t)\bigr) \bigr) + \mathsf P \bigl( \sigma< t,
X_\sigma\notin\bigl(c(t),b(t)\bigr) \bigr)\nonumber
\\
&\le&\mathsf P \bigl( X_\sigma\in\bigl(c(t),b(t)\bigr)
\bigr) + \mathsf P \bigl( \sigma< t, \tau\le\sigma, X_\sigma\notin
\bigl(c(t),b(t)\bigr) \bigr)\nonumber
\\
& =& \mathsf P \bigl( X_\tau\in\bigl(c(t),b(t)\bigr) \bigr) +
\mathsf P \bigl( \sigma< t, \tau= \sigma, X_\sigma\notin\bigl(c(t),b(t)
\bigr) \bigr)
\\
&\le&\mathsf P \bigl( \tau<t, X_\tau\in\bigl(c(t),b(t)\bigr)
\bigr) + \mathsf P \bigl( \tau< t, X_\tau\notin\bigl(c(t),b(t)\bigr)
\bigr)\nonumber
\\
& =& \mathsf P( \tau<t )\nonumber
\end{eqnarray}
for all $t>0$, proving the claim.
\end{pf}

%co12 #&#
\begin{cor}[(Uniform integrability)]\label{cor12} In the setting
of Theorem~\ref{teo1} let $\tau= \tau_{b,c}$ with $c = -\infty$ if
$\operatorname{supp}(\mu) \subseteq\mathbb{R}_+$ and $b = +\infty$
if $\operatorname{supp}(\mu)
\subseteq\mathbb{R}_-$.
%
%
%e4.5 #&#
%e4.6 #&#
%e4.7 #&#
\begin{eqnarray}
\label{45} &&\mbox{If $\int x \mu(dx) = 0$, then $\{B_{t \wedge\tau} \vert t \ge0 \}$ is uniformly integrable}.
\\
\label{46} &&\mbox{If $0 < \int x \mu(dx) < +\infty$, then $\{ B_{t \wedge\tau}^+ \vert t \ge0 \bigr\}$ is uniformly integrable}.
\\
\label{47} &&\mbox{If $-\infty< \int x \mu(dx) < 0$, then $\{ B_{t \wedge\tau}^- \vert t \ge0 \bigr\}$ is uniformly integrable}.
\end{eqnarray}
\end{cor}

\begin{pf}
Statement (\ref{45}) follows by combining
Proposition~\ref{prop11} above and Theorem 3 in \cite{Mo}, page 1294.
Statements (\ref{46}) and (\ref{47}) follow by combining
Proposition~\ref{prop11} above and Theorem 3 in \cite{CH}, page 397. This
completes the proof.
\end{pf}

%pr13 #&#
\begin{prop}[(Finiteness)]\label{prop13}
In the setting of
Theorem~\ref{teo1} suppose that $\operatorname{supp}(\mu) \cap\mathbb{R}_+
\ne\varnothing$ and
$\operatorname{supp}(\mu) \cap\mathbb{R}_- \ne\varnothing$.
%
%
%e4.8 #&#
%e4.9 #&#
\begin{eqnarray}
\qquad&&\mbox{If $\sup\operatorname{supp}(\mu) < \infty$, then there exists $T>0$ such that $b(t) = +\infty$}
\nonumber\\[-12pt]\label{48} \\[-12pt]
&&\mbox{for all $t>T$ if and only if $-\infty\le\int x \mu(dx) < 0$.}\nonumber
\\
&&\mbox{If $\inf\operatorname{supp}(\mu) > -\infty$, then there exists $T>0$ such that $c(t) = -\infty$}
\nonumber\\[-12pt]\label{49} \\[-12pt]
&&\mbox{for all $t>T$ if and only if $0 < \int x \mu(dx) \le+\infty$.}\nonumber
\end{eqnarray}
\end{prop}

\begin{pf}
It is enough to prove (\ref{48}) since (\ref{49})
then follows by symmetry. For this, suppose first that $b(t) =
+\infty$ for all $t>T$ with some minimal $T>0$. Since $\sup
\operatorname{supp}(\mu) < \infty$ we know that $b(T)< \infty$. Set
$b_1(t) =
b(t)$ for $t \in(0,T]$ and $b_1(t) = b(T)$ for $t>T$. Set $c_1(t) =
c(t)$ for $t \in(0,T]$ and $c_1(t) = c(T)$ for $t>T$ (recall that
$c$ must be finite valued). Then $\vert B_{t \wedge\tau_{b_1,c_1}}
\vert\le b(T) \vee(-c(T)) < \infty$ for all $t \ge0$ so that
$\{ B_{t \wedge\tau_{b_1,c_1}} \vert t \ge0 \}$ is
uniformly integrable and hence $\mathsf E B_{\tau_{b_1,c_1}} = 0$. Note
that $B_{\tau_{b,c}} \le B_{\tau_{b_1,c_1}}$ and moreover
$B_{\tau_{b,c}} < B_{\tau_{b_1,c_1}}$ on the set of a strictly
positive $\mathsf P$-measure where $B$ hits $b_1$ after $T$ before
hitting $c_1$. This implies that $\mathsf E B_{\tau_{b,c}} < \mathsf E
B_{\tau_{b_1,c_1}} = 0$ as claimed in (\ref{48}) above.

Conversely, suppose that $\mathsf E B_{\tau_{b,c}} < 0$, and consider first
the case when $c(t) = -\infty$ for $t>T$ with some $T>0$ at which
$c(T) > -\infty$. Set $c_1(t) = c(t)$ for $t \in(0,T]$ and $c_1(t)
= c(T)$ for $t>T$. Since $B_{\tau_{b,c_1}} \le\sup\operatorname
{supp}(\mu) <
\infty$ when $b$ is finite valued we see that $\vert B_{t \wedge
\tau_{b,c_1}} \vert\le\sup\operatorname{supp}(\mu) \vee(-c(T))
< \infty$
for all $t \ge0$ so that $\{ B_{t \wedge\tau_{b,c_1}} \vert
t \ge0 \}$ is uniformly integrable and hence $\mathsf E
B_{\tau_{b,c_1}} = 0$. Note that $B_{\tau_{b,c}} \ge
B_{\tau_{b,c_1}}$ so that $\mathsf E B_{\tau_{b,c}} \ge0$, and this
contradicts the hypothesis. Next consider the case when $c(t) >
-\infty$ for all $t \ge0$. Set $c_n(t) = c(t)$ for $t \in(0,n]$
and $c_n(t) = -\infty$ for $t>n$ with $n \ge1$. Set $d_n(t) = c(t)$
for $t \in(0,n]$ and $d_n(t) = c(n)$ for $t>n$ with $n \ge1$. Then
as above $\mathsf E B_{\tau_{b,d_n}} = 0$ and since $B_{\tau_{b,c_n}}
\ge
B_{\tau_{b,d_n}}$, it follows that $\mathsf E B_{\tau_{b,c_n}} \ge0$ for
all $n \ge1$. Moreover, since $B_{\tau_{b,c_n}} \le\sup
\operatorname{supp}(\mu) < \infty$ for all $n \ge1$ when $b$ is
finite valued by
Fatou's lemma, we get
%
%
%e4.10 #&#
\begin{equation}
\label{410} \mathsf E B_{\tau_{b,c}} = \mathsf E \lim_{n \rightarrow
\infty}
B_{\tau_{b,c_n}} \ge\limsup_{n \rightarrow\infty} \mathsf E B_{\tau_{b,c_n}}
\ge0,
\end{equation}
and this contradicts the hypothesis. Thus in both cases we see that
$b$ cannot be finite valued, and this completes the proof.
\end{pf}

%s5 #&#
\section{Optimality}\label{sec5}
%%%%%%%%%%%%%%%%%%%%

In this section we show that the stopping time from Theorem~\ref{teo1} has
the smallest truncated expectation among all stopping times that
embed $\mu$ into $B$. The same optimality result for stopping times
arising from the filling scheme when their means are finite was
derived by Chacon (\cite{Chh}, page 34), using a different method of
proof. The proof we present below is based on a recent proof of
Rost's optimality result \cite{Ro-2} given by Cox and Wang
\cite{CW}, Section~5. The verification technique we employ avoids
stochastic calculus and invokes a general martingale/Markovian
result to describe the supermartingale structure. This technique
applies in the setting of Corollary~\ref{cor8} as well and should also be of
interest in other/more general settings of this kind.

%th14 #&#
\begin{teo}\label{teo14}
In the setting of Theorem~\ref{teo1} or Corollary~\ref{cor8} let $\tau= \tau_{b,c}$ with $c = -\infty$ if $\operatorname
{supp}(\mu)
\subseteq\mathbb{R}_+$ and $b = +\infty$ if $\operatorname
{supp}(\mu) \subseteq\mathbb{R}_-$. If
$\sigma$ is any stopping time such that $B_\sigma\sim B_\tau$, then
we have
%
%
%e5.1 #&#
\begin{equation}
\label{51} \mathsf E (\tau\wedge T) \le\mathsf E (\sigma\wedge T)
\end{equation}
for all $T>0$.
\end{teo}

\begin{pf} Let $\mathsf P_{ t,x}$ denote the probability measure
under which $\mathsf P_{ t,x}(X_t=x) = 1$, and consider the function $H$
defined by
%
%
%e5.2 #&#
\begin{equation}
\label{52} H(t,x) = \mathsf P_{ t,x} ( \tau\le T )
\end{equation}
for $(t,x) \in[0,T] \times\mathbb{R}$ with $T>0$ given and fixed.
Extend $H$ outside $[0,T]$ by setting $H(t,x)=0$ for $t>T$ and $x
\in\mathbb{R}$. Define the (right) inverse $\rho$ of $b$ and $c$ by setting
\begin{eqnarray}
\label{53} \rho(x) & =& \inf\bigl\{ t>0 \vert b(t) \ge x \bigr\} \qquad
\mbox{if } x \ge b(0+)
\nonumber\\[-8pt]\\[-8pt]
\nonumber
& =& \inf\bigl\{ t>0 \vert c(t) \le x \bigr\} \qquad\mbox{if } x \le
c(0+).\nonumber
\end{eqnarray}
Then $x \mapsto\rho(x)$ is right-continuous and increasing on
$[b(0+),\infty)$ and left-continuous and decreasing on
$(-\infty,c(0+)]$. Set $D = (-\infty,c(0+)] \cup[b(0+),\infty)$ to
denote the domain of $\rho$, and note that $\rho(x) \ge0$ for all $x
\in D$.

\begin{longlist}[1.]
\item[1.] For $x \in D$ such that $\rho(x) \le T$ and $t \le\rho(x)$, we
have $H(s,x) = 1$ for all $s \in[t,\rho(x)]$. Hence we see that the
following identity holds
%
%
%e5.3 #&#
\begin{equation}
\label{54} \rho(x) - t = \int_t^{\rho(x)} H(s,x)
\,ds
\end{equation}
whenever $t \le\rho(x) \le T$. Since $H \le1$, we see that this
identity extends as
%
%
%e5.4 #&#
\begin{equation}
\label{55} \rho(x) - t \le\int_t^{\rho(x)}
H(s,x) \,ds
\end{equation}
for $\rho(x) < t \le T$. Since $\rho(x) - t = (T - t)^+ - (T -
\rho(x))^+$ for $t \vee\rho(x) \le T$ and $H(s,x) = 0$ for $s > T$,
it is easily verified using the same arguments as above that
(\ref{54}) and (\ref{55}) yield
%
%
%e5.5 #&#
\begin{equation}
\label{56} (T - t)^+ \le\int_t^{\rho(x) \wedge T} H(s,x)
\,ds + \bigl(T - \rho(x)\bigr)^+
\end{equation}
for all $t \ge0$ and $x \in D$. Let us further rewrite (\ref{56})
as follows:
%
%
%e5.6 #&#
\begin{equation}
\label{57} (T - t)^+ \le F(t,x) + G(x),
\end{equation}
where the functions $F$ and $G$ are defined by
%
%
%e5.7 #&#
%e5.8 #&#
\begin{eqnarray}
\label{58} F(t,x) &=& \int_t^T H(s,x) \,ds,
\\
\label{59}  G(x) &=& \bigl(T - \rho(x)\bigr)^+ - \int_{\rho(x) \wedge T}^T
H(s,x) \,ds
\end{eqnarray}
for $t \ge0$ and $x \in D$.
\end{longlist}

\begin{longlist}[2.]
\item[2.] It is easily seen from definitions of $\tau$ and $\rho$ (using
that $b$ and $c$ are increasing and decreasing, resp.) that
$\rho(X_\tau) \ge\tau$. Combining this with the fact that \mbox{$H(s,x) =
1$} for all $s \in[t,\rho(x) \wedge T]$ and $x \in D$, we see that
equality in (\ref{56}) is attained at $(\tau,X_\tau)$. Since
(\ref{57}) is equivalent to (\ref{56}), it follows that
%
%
%e5.9 #&#
\begin{equation}
\label{510} (T - \tau)^+ = F(\tau,X_\tau) + G(X_\tau).
\end{equation}
We now turn to examining (\ref{57}) for other stopping times.
\end{longlist}

\begin{longlist}[3.]
\item[3.] To understand the structure of the function $F$ from (\ref{58}),
define
%
%
%e5.10 #&#
\begin{equation}
\label{511} D_t = \bigl\{ (s,x) \in\mathbb{R}_+ \times
\mathbb{R} \vert x \ge b(t + s) \mbox{ or } x \le c(t + s) \bigr
\},
\end{equation}
and note by time-homogeneity of $X$ that
%
%
%e5.11 #&#
\begin{equation}
\label{512} H(t,x) = \mathsf P_{ t,x} ( \tau\le T ) = \mathsf
P_{ x} ( \tau_t \le T - t )
\end{equation}
for $(t,x) \in[0,T] \times\mathbb{R}$ where we set
%
%
%e5.12 #&#
\begin{equation}
\label{513} \tau_t = \inf\{ s>0 \vert X_s \in
D_{t+s} \}
\end{equation}
with respect to the probability measure $\mathsf P_{ x}$ under which
$\mathsf P_{ x}(X_0=x)=1$. Hence we see that
%
%
%e5.13 #&#
\begin{eqnarray}
\label{514} F(t,x) & =& \int_t^T H(s,x) \,ds =
\int_t^T \mathsf P_{ x} (
\tau_s \le T - s ) \,ds
\nonumber\\[-8pt]\\[-8pt]
& =& \int_0^{T-t} \mathsf
P_{ x} ( \tau_{T-s} \le s ) \,ds = \mathsf E _x
\int_0^{T-t} Z_s \,ds\nonumber
\end{eqnarray}
for $(t,x) \in[0,T] \times\mathbb{R}$ where we set
%
%
%e5.14 #&#
\begin{equation}
\label{515} Z_s = I ( \tau_{T-s} \le s )
\end{equation}
for $s \in[0,T - t]$. Noting that each $Z_s$ is
${\mathcal F}_s$-measurable where ${\mathcal F}_s = \sigma(X_r \vert
0 \le r \le
s)$, we can now invoke a general martingale/Markovian result and
conclude that
%
%
%e5.15 #&#
\begin{equation}
\label{516} M_t:= F(t,X_t) + \int
_0^t Z_s \,ds
\end{equation}
is a martingale with respect to ${\mathcal F}_t$ for $t \in[0,T]$. Indeed,
for this note that by the Markov property of $X$, we have
%
%
%e5.16 #&#
\begin{eqnarray}
\label{517} \mathsf E _x(M_{t+h} |{\mathcal
F}_t) & =& \mathsf E _x \biggl( F(t + h,X_{t+h})
+ \int_0^{t+h} Z_s \,ds \Big|{\mathcal
F}_t \biggr)\nonumber
\\
& =& \mathsf E _x \biggl( \mathsf E _{X_{t+h}} \biggl(
\int_0^{T-t-h} Z_s \,ds \biggr) + \int
_0^{t+h} Z_s \,ds \Big|{\mathcal
F}_t \biggr)\nonumber
\\
& =& \mathsf E _x \biggl( \mathsf E _x \biggl(
\int_0^{T-t-h} Z_s \,ds \circ
\theta_{t+h} \Big|{\mathcal F}_{t+h} \biggr) + \int
_0^{t+h} Z_s \,ds \Big|{\mathcal
F}_t \biggr)\nonumber
\\
& =& \mathsf E _x \biggl( \int_0^T
Z_s \,ds \Big|{\mathcal F}_t \biggr) = \mathsf E
_x \biggl( \int_t^T
Z_s \,ds \Big|{\mathcal F}_t \biggr) + \int
_0^t Z_s \,ds
\\
& =& \mathsf E _x \biggl( \int_0^{T-t}
Z_s \,ds \circ\theta_t \Big|{\mathcal F}_t
\biggr) + \int_0^t Z_s \,ds\nonumber
\\
& =& \mathsf E _{X_t} \biggl( \int_0^{T-t}
Z_s \,ds \biggr) + \int_0^t
Z_s \,ds\nonumber
\\
& =& F(t,X_t) + \int_0^t
Z_s \,ds = M_t\nonumber
\end{eqnarray}
for all $0 \le t \le t + h \le T$, showing that (\ref{516}) holds
as claimed. Extend the martingale $M$ to $(T,\infty)$ by setting
$M_t=M_T$ for $t>T$. Since $F(t,x)=0$ for $t>T$ and $x \in\mathbb{R}$, this
is equivalent to setting $Z_s=0$ for $s>T$ in (\ref{516}) above.
Since $Z_s \ge0$ for all $s \ge0$ we see from (\ref{516}) that
$F(t,X_t)$ is a supermartingale with respect to ${\mathcal F}_t$ for $t
\ge
0$.
\end{longlist}

\begin{longlist}[4.]
\item[4.] We next note that
%
%
%e5.17 #&#
\begin{equation}
\label{518} \int_0^{t \wedge\tau} Z_s \,ds =
0
\end{equation}
for all $t \ge0$. Indeed, this is due to the fact that $\tau_{T-s}
= \inf\{ r>0 | X_r \in D_{T-s+r} \} \ge\inf\{
r>0 | X_r \in D_0 \} = \tau$ for all $s \in[0,\tau\wedge
T)$ since $b$ is increasing and $c$ is decreasing. Hence from
(\ref{515}) we see that $Z_s = 0$ for all $s \in[0,\tau)$, and this
implies~(\ref{518}) as claimed. Combining (\ref{516}) and
(\ref{518}) we see that $F(t \wedge\tau,X_{t \wedge\tau})$ is a
martingale with respect to ${\mathcal F}_{t \wedge\tau}$ for $t \ge0$.
\end{longlist}

\begin{longlist}[5.]
\item[5.] Taking now any stopping time $\sigma$ such that $X_\sigma\sim
X_\tau$ it follows by (\ref{510}), (\ref{518}), (\ref{516}) and
(\ref{57}) using the optional sampling theorem that
%
%
%e5.18 #&#
\begin{eqnarray}
\label{519} \qquad \mathsf E (T - \tau)^+ & =& \mathsf E F(\tau,X_\tau) +
\mathsf E G(X_\tau) = \mathsf E M_\tau+ \mathsf E
G(X_\sigma)
\nonumber\\[-8pt]\\[-8pt]
& =& \mathsf E M_\sigma+ \mathsf E G(X_\sigma) \ge
\mathsf E F(\sigma,X_\sigma) + \mathsf E G(X_\sigma) \ge\mathsf
E (T - \sigma)^+.\nonumber
\end{eqnarray}
Noting that $\mathsf E (T - \tau)^+ = T - \mathsf E (\tau\wedge T)$
and $\mathsf E (T
- \sigma)^+ = T - \mathsf E (\sigma\wedge T)$, we see that this is
equivalent to (\ref{51}), and the proof is complete.\quad\qed
\end{longlist}\noqed
\end{pf}

%re15 #&#
\begin{rem}\label{rem15}
In the setting of Theorem~\ref{teo1} if $\int x^2
\mu(dx) < \infty$, then $\mathsf E B_\tau^2 < \infty$ and hence
$\mathsf E \tau<
\infty$ since $\tau$ is minimal (Section~\ref{sec4}). If moreover $\mathsf E
\sigma
< \infty$, then by It\^o's formula and the optional sampling theorem,
we know that $\mathsf E \sigma= \mathsf E \tau$. When $\int x^2 \mu
(dx) =
\infty$, however, it is not clear a priori whether the ``expected
waiting time'' for $\tau$ compares favourably with the ``expected
waiting time'' for any other stopping time $\sigma$ that embeds $\mu$
into $B$. The result of Theorem~\ref{teo14} states the remarkable fact that
$\tau$ has the smallest truncated expectation among all stopping
times $\sigma$ that embed $\mu$ into $B$ (note that this fact is
nontrivial even when $\mathsf E \tau$ and $\mathsf E \sigma$ are
finite). It is
equally remarkable that this holds for all laws $\mu$ with no extra
conditions imposed.

The optimality result of Theorem~\ref{teo14} extends to more general concave
functions using standard techniques.
\end{rem}

%co16 #&#
\begin{cor}[(Optimality)]\label{cor16}
In the setting of Theorem~\ref{teo1} or Corollary~\ref{cor8}, let $\tau= \tau_{b,c}$ with $c = -\infty$ if
$\operatorname{supp}(\mu) \subseteq\mathbb{R}_+$ and $b = +\infty$
if $\operatorname{supp}(\mu)
\subseteq\mathbb{R}_-$, and let $F\dvtx \mathbb{R}_+ \rightarrow
\overline{\mathbb{R}}$ be a concave
function such that $\mathsf E F(\tau)$ exists. Then we have
%
%
%e5.19 #&#
\begin{equation}
\label{520} \mathsf E F(\tau) \le\mathsf E F(\sigma)
\end{equation}
for any stopping time $\sigma$ such that $X_\sigma\sim X_\tau$.
\end{cor}

\begin{pf}
By (\ref{51}) we know that
%
%
%e5.20 #&#
\begin{equation}
\label{521} \int_0^t \mathsf P(\tau> s)
\,ds \le\int_0^t \mathsf P(\sigma> s) \,ds
\end{equation}
for all $t \ge0$. It is easy to check using Fubini's theorem that
for any nonnegative random variable $\rho$ we have
%
%
%e5.21 #&#
\begin{equation}
\label{522} \mathsf E F(\rho) = F(0) - \int_0^\infty\!\!
\int_0^t \mathsf P(\rho> s) \,ds
F'(dt)
\end{equation}
whenever $F$ is a concave function satisfying $t F'(t) \rightarrow
0$ as $t \downarrow0$ and $F'(t) \rightarrow0$ as $t \rightarrow
\infty$ where $F'$ denotes the right derivative of $F$. Applying
(\ref{522}) to $\tau$ and $\sigma$, respectively, recalling that
$F'(dt)$ defines a negative measure, and using (\ref{521}) we get
(\ref{520}) for those functions $F$. The general case then follows
easily by tangent approximation (from the left) and/or truncation
(from the right) using monotone convergence.
\end{pf}

%re17 #&#
\begin{rem}\label{rem17}
In addition to the temporal optimality of $b$
and $c$ established in~(\ref{520}), there also exists their spatial
optimality arising from the optimal stopping problem
%
%
%e5.22 #&#
\begin{equation}
\label{523} \sup_{0 \le\tau\le T} \mathsf E \biggl( |
B_\tau|- 2 \int_0^{B_\tau}
F_\mu(x) \,dx \biggr),
\end{equation}
where $F_\mu$ denotes the distribution function of $\mu$. Indeed
McConnell (\cite{Mc}, Section~5), shows that (under his conditions) the
optimal stopping time in (\ref{523}) equals
%
%
%e5.23 #&#
\begin{equation}
\label{524} \tau_* = \inf\bigl\{ t \in[0,T] | B_t \ge b(T -
t) \mbox{ or } B_t \le c(T - t) \bigr\},
\end{equation}
where $b$ and $c$ are functions from Theorem~\ref{teo1} (compare (\ref{523})
with the optimal stopping problem derived in \cite{Pe}). This can be
checked by the It\^o--Tanaka formula and the optional sampling theorem
from the local time reformulation of (\ref{523}) that reads
%
%
%e5.24 #&#
\begin{equation}
\label{525} \sup_{0 \le\tau\le T} \mathsf E \biggl( \int
_{\mathbb{R}} \ell_\tau^x \nu(dx) - \int
_{\mathbb{R}} \ell_\tau^x \mu(dx) \biggr),
\end{equation}
where $\ell$ is the local time of $B$, and $\nu$ is a probability
measure on $\mathbb{R}$ such that $\operatorname{supp}(\nu)
\subseteq[-p,q]$ with
$\mu([-p,q])=0$ for some $p>0$ and $q>0$. Since the existence and
uniqueness result of Theorems~\ref{teo1}~and~\ref{teo10} with $B_0 \sim\nu$
remain valid in this case as well (recall Remark~\ref{rem2} and the beginning
of Section~\ref{sec3}), we see that McConnell (\cite{Mc}, Section~5), implies
that (under his conditions) the resulting stopping time (\ref{524})
is optimal in (\ref{525}).
\end{rem}

%\begin{appendix}
%\section{}
%\end{appendix}

% zodis "Acknowledgments" paliekamas pagal autoriu
%\section*{Acknowledgments}

%\begin{supplement}[id=suppA]
%\sname{Supplement A}
%\stitle{}
%\slink[doi]{10.1214/00-AOPXXXXSUPP} %[doi,text={...}] - jei reikia
%suskaldyti doi
%\sdatatype{.pdf}
%\sfilename{aopXXXX\_supp.pdf}
%\sdescription{}
%\end{supplement}

% imsref loaded by linak, 2014-06-17 16:01:07
%

\printaddresses
\end{document}